\documentclass{amsart}

\usepackage[utf8]{inputenc}
\usepackage{amssymb}
\usepackage{amsmath}
\usepackage{amstext}
\usepackage{amsthm}
\usepackage{enumitem}
\usepackage[colorlinks=true]{hyperref}
\usepackage{xcolor}
\definecolor{lightPKU}{cmyk}{0,1,1,0.2}

\usepackage{graphicx}
\usepackage{tabularx}
\usepackage{booktabs}
\usepackage{float}
\usepackage{bibref}
\usepackage{tikz}
\usepackage{tikz-cd}
\tikzcdset{scale cd/.style={every label/.append style={scale=#1},
		cells={nodes={scale=#1}}}}

\usepackage[linesnumbered,ruled,vlined]{algorithm2e}

\SetCommentSty{mycommfont}

\usepackage{mathrsfs}
\usepackage{complexity}
\usepackage{subcaption}

\newcommand{\tr}{\mathop{\mathrm{tr}}}
\newcommand{\rank}{\mathop{\mathrm{rank}}}
\newcommand{\spec}{\mathop{\mathrm{Spec}}}
\newcommand{\Frac}{\mathop{\mathrm{Frac}}}
\newcommand{\sign}{\mathop{\mathrm{sign}}}
\newcommand{\ideal}[1]{\langle#1\rangle}
\newcommand{\bs}[1]{\boldsymbol{#1}}
\newcommand{\lt}{{\mathtt{LT}}}
\newcommand{\lm}{{\mathtt{LM}}}
\newcommand{\lc}{{\mathtt{LC}}}
\newcommand{\res}{{\mathtt{Res}}}

\makeatletter
\newcommand\xleftrightarrow[2][]{%
	\ext@arrow 9999{\longleftrightarrowfill@}{#1}{#2}}
\newcommand\longleftrightarrowfill@{%
	\arrowfill@\leftarrow\relbar\rightarrow}
\makeatother

\newcounter{algoline}

\AtBeginEnvironment{algorithm}{\setcounter{algoline}{0}}

\newtheorem{theorem}{Theorem}[section]
\newtheorem{lemma}[theorem]{Lemma}

\newtheorem{corollary}[theorem]{Corollary}
\theoremstyle{definition}
\newtheorem{definition}[theorem]{Definition}
\newtheorem{example}{Example}

\newtheorem*{problem}{Problem}
\theoremstyle{remark}
\newtheorem{remark}{\indent Remark}
\newtheorem*{remark*}{Remark}

\begin{document}

	\title{A geometric approach to cylindrical algebraic decomposition}
	\author{Rizeng Chen}
	
	\address{School of Mathematical Sciences, Peking University, 100871, Beijing, China}
	\email{xiaxueqaq@stu.pku.edu.cn}
\begin{abstract}
	 Cylindrical algebraic decomposition is a classical construction in real algebraic geometry. Although there are many algorithms to compute a cylindrical algebraic decomposition, their practical performance is still very limited. In this paper, we revisit this problem from a more geometric perspective, where the construction of cylindrical algebraic decomposition is related to the study of morphisms between real varieties. It is showed that the geometric fiber cardinality (geometric property) decides the existence of semi-algebraic continuous sections (semi-algebraic property). As a result, all equations can be systematically exploited in the projection phase, leading to a new simple algorithm whose efficiency is demonstrated by experimental results.
\end{abstract}
\thanks{This work is supported by National Key R\&D Program of China (No. 2022YFA1005102).}
\maketitle

\section{Introduction}
	\noindent\textbf{Background.}  Cylindrical algebraic decomposition is a fundamental tool in real algebraic geometry. Roughly speaking, the goal of it is to decompose semi-algebraic sets into cylindrically arranged cells that are homeomorphic to open balls of various dimensions. The study of cylindrical algebraic decomposition began with Collins's landmark paper \cite{collins1975quantifier}. Collins's original motivation was to develop an algorithm to efficiently solve the quantifier elimination problem over the reals, since Tarski's method \cite{tarski1951decision} is infeasible even for very small instances. But later it is proved that the underlying idea also plays an important role in studying the geometry of semi-algebraic sets. Collins's strategy for cylindrical algebraic decomposition consists of two phases: repeated projection and lifting. Often the projection phase takes a long time, leading to a bottleneck. Hence there have been exhaustive efforts on improving the projection phase, e.g.\ \cite{hong1990improvement}, \cite{lazard1994improved}, \cite{mccallum1998improved} and \cite{brown2001improved}.
	\bigskip
	
	\noindent\textbf{Challenge.} It was soon noticed that some information can be used to accelerate the computation of c.a.d., especially the equations. For some applications, there are already equations in the input (e.g.\ Real Root Classification, Automatic Theorem Proving). Also, equations are generated during the projection phase. 
		
\begin{figure}[hbtp]
		\centering
		\begin{subfigure}[t]{0.35\textwidth}
		\centering
			\includegraphics[width=\linewidth]{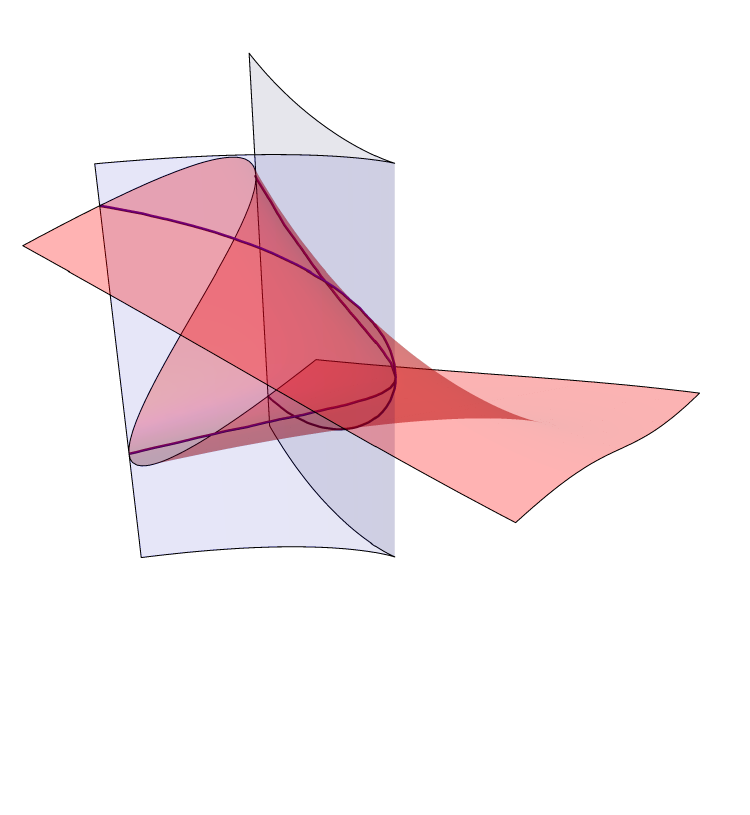}
			\caption{The curve as intersection of two surfaces} 
			\label{fig:depressed-cubic-over-discriminant}
		\end{subfigure}
		~
		\begin{subfigure}[t]{0.35\textwidth}
			\includegraphics[width=\linewidth]{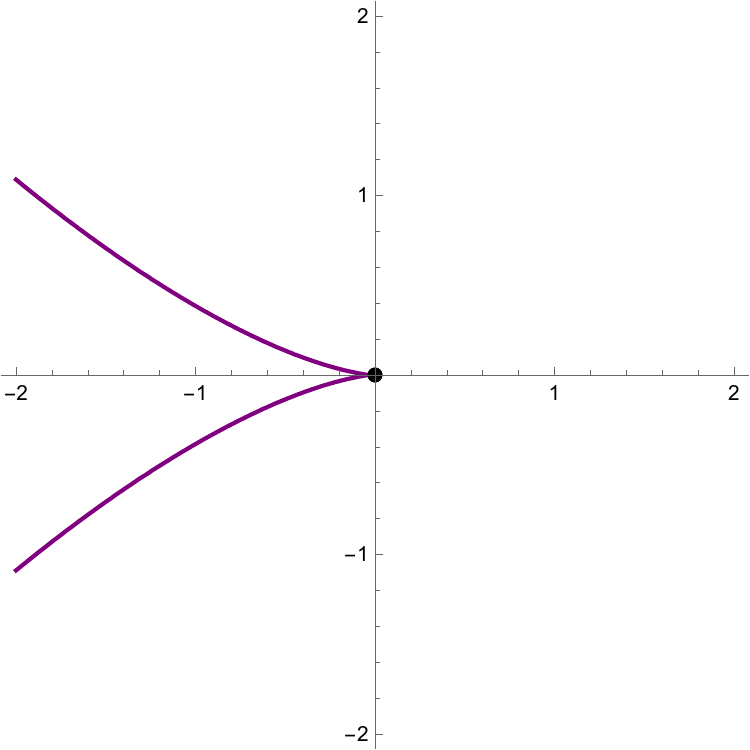}
			\caption{The projection image}
			\label{fig:discriminant-origin}
		\end{subfigure}
		~
		\caption{Equations naturally arise}
	\end{figure}

		To illustrate this observation, let us consider a toy polynomial system $x^3+px+q=4p^3+27q^2=0$ in parameters $p,q$ and variable $x$. 
		It defines a curve in the three-dimensional space with two irreducible components intersecting at the origin. The curve is depicted in Figure \ref{fig:depressed-cubic-over-discriminant} in purple, as the intersection of two surfaces defined by $x^3+px+q=0$ in red and $4p^3+27q^2=0$ in blue, respectively. One wants to study the number of real roots of the system. This is a classical application of cylindrical algebraic decomposition. Certainly there are equations in the input and extra equations are generated during the projection ($p=q=0$, reflecting the intersection of two irreducible components). See Figure \ref{fig:discriminant-origin}.
				
	For larger non-trivial inputs, this phenomena become much more frequent: many equations are already in the input and many additional equations are systematically generated during the projection phase. Therefore it is an important challenge to exploit these equations.
	
	\bigskip
	\noindent\textbf{Previous Results.} There have been intensive efforts regarding the issue and a significant process is made to utilize the equational constraints in the input in \cite{mccallum1999projection,mccallum2001propagation}, based on McCallum's projection operator \cite{mccallum1998improved}. The main idea is to select an equation (if there are any) at each stage of projection and decompose the hypersurface defined by the equation only and implicit equations are discovered by resultants of known equations. Later McCallum and Brown generalized this to two equations \cite{brown2005using,mccallum2009delineability}. It is also found that savings can be made during the lifting phase \cite{england2015improving}. For a recent summary on utilizing Equational Constraints in CAD, we refer readers to \cite{england2020cylindrical}. 

	However, there is still more to be done. 
	\begin{enumerate}
		\item The previous results did not make full use of the equations. Only one or two equations are used as ``pivot'' in each stage of projection. And they ignored many equations emerging in the middle of projection (e.g.\ the origin in Figure \ref{fig:discriminant-origin}).
		
		To further demonstrate this, let us consider another toy equation $px+q=0$, where $p,q$ are parameters and $x$ is the unknown. Its locus is depicted in Figure \ref{fig:linear-eqn}. Clearly, the equation has exactly one solution when $p\ne 0$, zero solution when $p=0$ but $q\ne 0$, and infinite solutions when $p=q=0$.
		Therefore during the projection, equations $p=0$ and $p=q=0$ naturally arise. However, the latter is neglected in classical algorithms. Instead, they compute a ${p,q}$-sign-invariant c.a.d.\ of the $pq$-plane. As a result, the $pq$-plane is divided into 9 cylindrical cells (four quadrants, four half-axes and the origin), although 5 is enough ($p>0$, $p<0$, the positive and negative half-axes of $q$, and the origin).
		\begin{figure}[htbp]
			\centering
			\begin{subfigure}[t]{0.35\textwidth}
				\includegraphics[width=\linewidth]{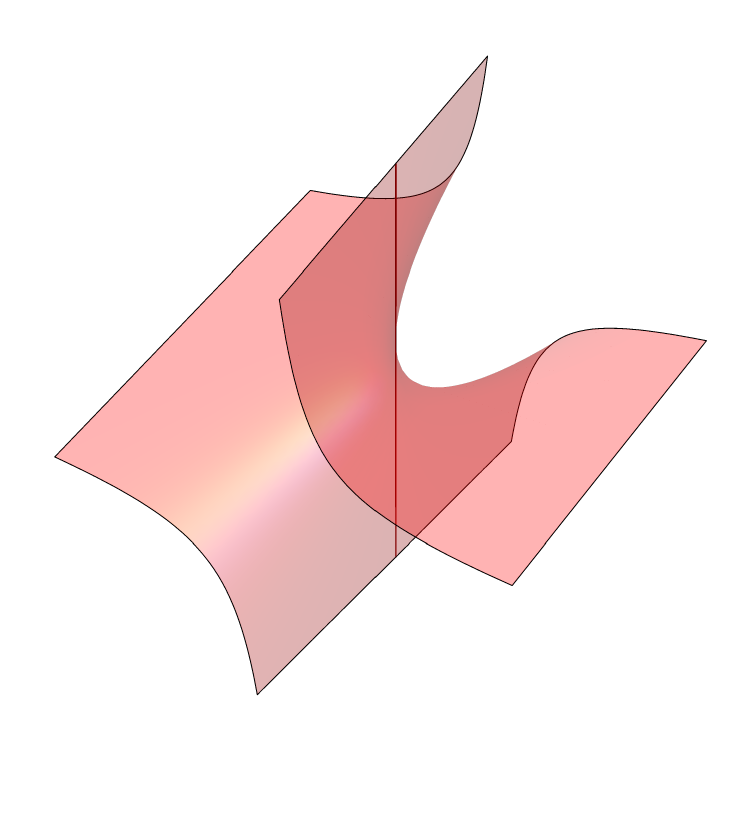}
				\caption{The surface $V(px+q)$, note that the $x$-axis is contained in the surface}
				\label{fig:linear-eqn}
			\end{subfigure}
			~
			\begin{subfigure}[t]{0.35\textwidth}
				\includegraphics[width=\linewidth]{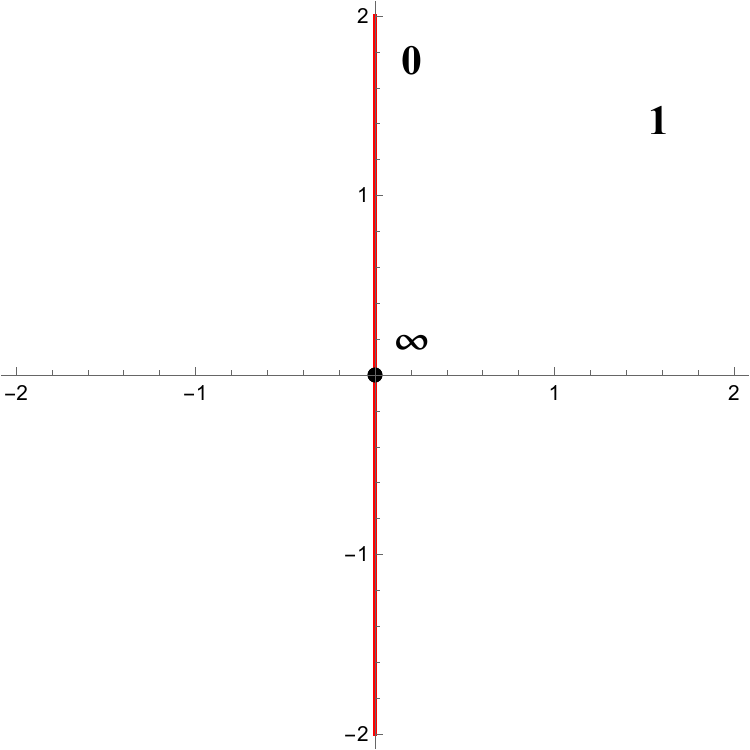}
				\caption{Number of solutions in $x$}
				\label{fig:y-axis-origin}
			\end{subfigure}
			~
			\caption{Some equations are ignored}
		\end{figure}

		\item The result in \cite{mccallum2001propagation} works with irreducible polynomials. A complicated preprocessing needs to be done to apply their result.
		\item McCallum's projection operator \cite{mccallum1998improved} could fail in rare cases. A key property  ``order-invariance'' is not maintained during the lifting phase, if the projection of a hypersurface is not quasi-finite. As a result, the correctness is not guaranteed in this case. This is called the ``Nullification Problem'' in the literature. For such an example, we refer readers to Example \ref{ex-5-var}.
	\end{enumerate}

	%
	\bigskip
	\noindent\textbf{Contribution.} In this paper, we develop a new simple algorithm that exploits all the equations, based on the study of finite free morphisms of affine real varieties. 
	\begin{enumerate}
		\item All the equations (from both the input and the intermediate stages) are used.
		\item It does not rely on radical (square-free) computation nor primary decomposition (factorization), so it is conceptually simple.
		\item It is a complete algorithm that never fail (no nullification problem).
		\item Most importantly, it is generally faster than the existing methods. A benchmark is reported in Section \ref{subsect-experiments}.
	\end{enumerate}

	The key result enabling all these advances is Theorem \ref{thm-freeness-implies-geometric-delineability}, which shows that a finite free morphism $f:X\to Y$ of affine real varieties admits semi-algebraic continuous sections on regions having a constant geometric fiber size, if $X$ is a closed subscheme of the cylinder $\mathbb{A}_Y^1$. This is a non-trivial generalization of \cite[Theorem 1]{collins1975quantifier} in the algebro-geometric context, connecting a geometric invariant (cardinality of geometric fiber) and a semi-algebraic property (existence of semi-algebraic continuous sections). Thus the theorem itself provides a link between real and complex algebraic geometry.	
	
	The result is turned into an algorithm by an effective version of Grothendieck's Generic Freeness and a geometric fiber counting technique (Hermite Quadratic Form). Very roughly speaking, our algorithm starts with some varieties. Then we stratify their projection images into ``flat'' (fibers changes continuously) and ``unramified'' (fibers do not break into different branches) pieces. 
	
		We illustrate the above idea with a c.a.d.\ of the plane. The same idea applies to high-dimensional varieties. Let $f=x y^3 - y + x^2$, the goal is to compute an $f$-sign-invariant c.a.d.\ of the plane.
		\begin{figure}[htbp]
			\centering
			\begin{subfigure}[t]{0.3\textwidth}
				\includegraphics[width=\linewidth]{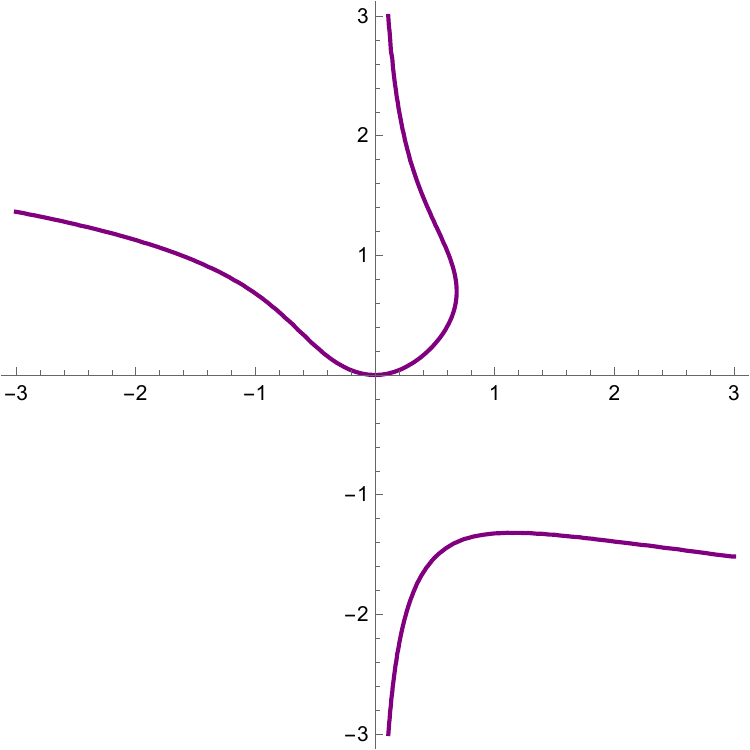}
				\caption{Locus of $f$} 
				\label{fig:2d-curve}
			\end{subfigure}
			~
			\begin{subfigure}[t]{0.3\textwidth}
				\includegraphics[width=\linewidth]{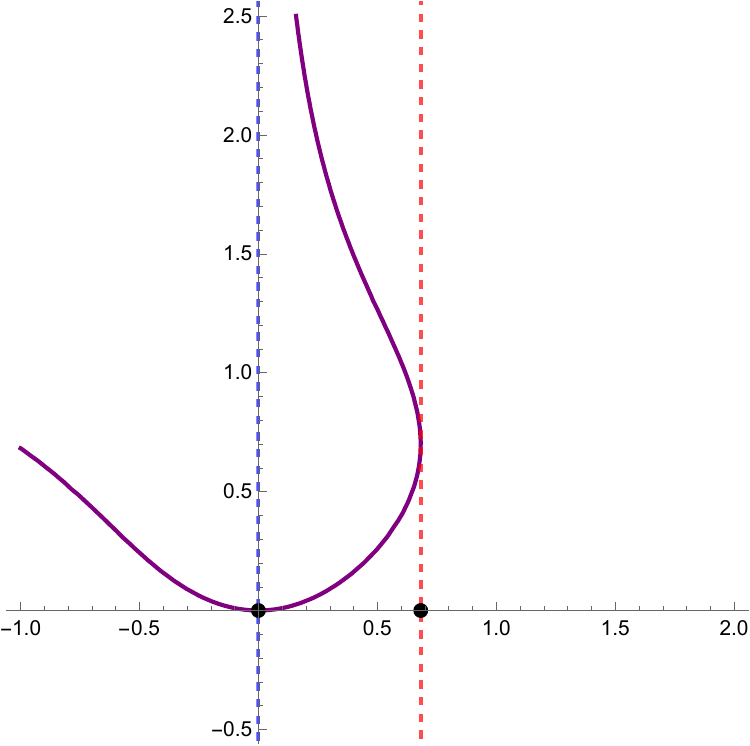}
				\caption{Points where ``flatness'' or ``unramifiedness'' fails} 
				\label{fig:2d-curve-nonflat-ramified}
			\end{subfigure}
			~			
			\begin{subfigure}[t]{0.3\textwidth}
				\includegraphics[width=\linewidth]{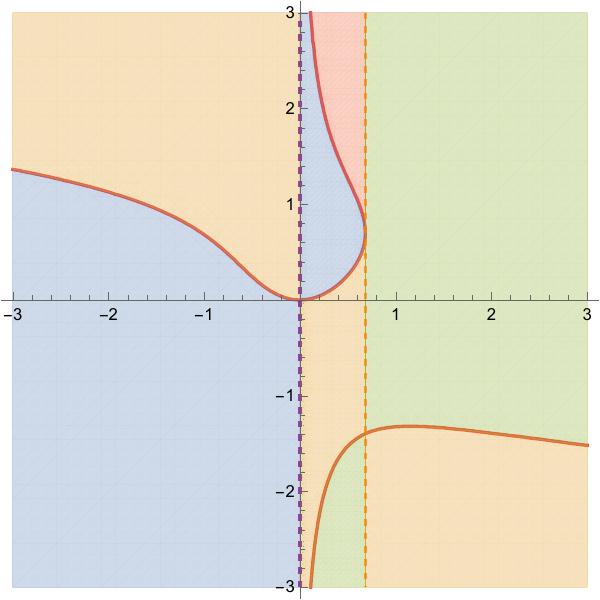}
				\caption{The final c.a.d.} 
				\label{fig:2d-cad}
			\end{subfigure}
			~
			\caption{An example of c.a.d.}
		\end{figure}
	
		We study the projection of the curve $V(f)$. It can be seen that the fibers change continuously when $x\ne 0$, so the affine line is split into two parts where the projection is ``flat''. Notice that the fiber ramifies when $27x^5-4=0$ (e.g.\ $x=\sqrt[5]{4/27}$ on the real axis). Therefore the $x\ne 0$ part is further divided according to whether $27x^5-4$ is zero or not.
		
		So the projection is ``flat'' and ``unramified'' over
		\begin{itemize}
			\item  $x\ne 0\wedge 27x^5-4\ne 0$, 
			\item $x\ne 0 \wedge 27x^5-4=0$, 
			\item $x=0$.
		\end{itemize}
		By Theorem \ref{thm-freeness-implies-geometric-delineability}, the projection of $V(f)$ admits sections over each semi-algebraically connected components of these regions. As a result, a c.a.d.\ is constructed over them, see Figure \ref{fig:2d-cad}.

	The crucial difference of our approach and the existing methods is the systematic use of equations, which requires a deeper understanding of morphisms of real varieties instead of just hypersurfaces. The algebro-geometric result Theorem \ref{thm-freeness-implies-geometric-delineability} is the essential ingredient to the new algorithm.
	
	\bigskip
	\noindent\textbf{Assumptions:} We assume the readers are familiar with the following basic concepts, which can be easily found in textbooks and literature:
	\begin{itemize}
		\item ordered fields, real closed fields, real closure, semi-algebraic set, semi-algebraic functions, Euclidean topology and semi-algebraic connectedness \cite{basu2008algorithms,bochnak2013real,mangolte2020real,scheiderer2024course} and
		\item monomial order $T$ (e.g.\ $\mathtt{grevlex}$ and block order), leading coefficient $\lc_T(f)$, leading term $\lt_T(f)$, leading monomial $\lm_T(f)$ and Gr\"{o}bner basis \cite{buchberger1965algorithmus,cox2013ideals}.
	\end{itemize} 
	\bigskip
	\noindent\textbf{Structure of the paper.} 
	Section \ref{sect-prob-stat} states the problem precisely.
	Section \ref{sect-preliminary} provides a quick review of affine scheme since it will play a crucial role in the subsequent theoretical development.
	In Section \ref{sect-fib-cls-of-gen-finite-mor}, a trace pairing technique (Hermite Quadratic Form) is generalized to count the geometric fiber of a finite free morphism of varieties. Then in Section \ref{sect-geometric-cad}, we build the geometric theory needed to study cylindrical algebraic decomposition and describe a new algorithm for cylindrical algebraic decomposition. Section \ref{sect-analysis} first compares our new result to existing methods, then reports on benchmarks and finally concludes the paper with a discussion on future works. 

\section{Problem Statement}	\label{sect-prob-stat}
In this section, we give the precise statement of the problem. For this, we will first present here the definition of cylindrical algebraic decomposition, which is essentially the same as \cite[Definition 5.1]{basu2008algorithms}. A minor difference is that here a c.a.d. does not necessarily cover the whole affine space $R^n$ (a natural generalization for applications like real root classification).
\begin{definition}
	Let $R$ be a real closed field.
	A \textbf{cylindrical algebraic decomposition} (c.a.d.) $\mathfrak{D}$ in $R^n$, is a sequence $S_1,\ldots,S_n$, where for each $1\leq i\leq n$:
	\begin{enumerate}
		\item The set $S_i$ consists of disjoint semi-algebraic sets in $R^i$, which are called the \textbf{cells} of level $i$.
		\item Each cell of level $1$ is either a point or an interval.
		\item For every $1< i\leq n$, and every $S\in S_i$, the \textbf{silhouette}
		$$\pi_{i-1}(S)=\{(y_1,\ldots,y_{i-1})\in R^{i-1}|\exists x\in R \text{ s.t. }(y_1,\ldots,y_{i-1},x)\in S\}$$ is a cell of level $i-1$.
		\item For every $1\leq i<n$, and every $S\in S_i$, there are finitely many semi-algebraic continuous functions $\xi_{S,1},\ldots,\xi_{S,l_S}:S\to R$ (the \textbf{sections}) such that any cell in $S_{i+1}$ lying over $S$ is either the graph of some $\xi_{S,j}$ 
		$$\left\{(y,\xi_{S,j}(y))\in S\times R\middle|y\in S\right\},$$
		or the \textbf{band} 
		$$\left\{(y,x)\in S\times R\middle|y\in S,\ \xi_{S,j}(y)<x<\xi_{S,j+1}(y)\right\},$$
		where we take $\xi_{S,0}(x)=-\infty$ and $\xi_{S,l_S+1}(x)=+\infty$.
	\end{enumerate}
	
	We say $\mathfrak{D}$ is \textbf{$\mathcal{F}$-sign-invariant} for a finite family of $n$-variate polynomials $\mathcal{F}$, if for every $f\in \mathcal{F}$ and every cell $\mathcal{C}$ of level $n$, $f$ is sign-invariant on $\mathcal{C}$.
	
	We say $\mathfrak{D}$ is \textbf{adapted to} a finite family of semi-algebraic sets $\{L_1,\ldots,L_s\}$ in $R^n$, if each $L_i$ is a union of level $n$ cells. Notice that the disjointness of cells implies that for each cell $S\in S_n$ of level $n$ and each semi-algebraic set $L_i$, either $S\subseteq L_i$ or $S\cap L_i=\varnothing$.		
	
	
\end{definition}

\begin{problem}
	Design an algorithm that
	
	\textbf{Inputs:} A list of varieties $V(I_1),\ldots,V(I_s)$ and
	
	\textbf{Outputs:} A cylindrical algebraic decomposition adapted to $\{V(I_1),\ldots,V(I_s)\}$.
\end{problem}
Note that a c.a.d.\ adapted to $\{V(0),V(f_1),\ldots,V(f_s)\}$ is an $\{f_1,\ldots,f_s\}$-sign-invariant c.a.d.. So the classical problem of computing a sign-invariant c.a.d.\ is a special case here.

Note also that it is easy to convert such an algorithm to accept semi-algebraic sets input.
\section{A Quick Review on Affine Schemes}\label{sect-preliminary}
Now we will begin to describe our approach tackling the problem. It will crucially use some fundamental geometric tools on Grothendieck's language of schemes. While it is true that most of our results can be formulated in the classical language of varieties, they will become very complicated and less intuitive. Since we are doing algebraic geometry over a non-algebraically closed field ($\mathbb{R}$ for instance) and prime/radical ideals are rare in computation, schemes should serve as the right settings (we do not require varieties to be integral for the same reason). 

Here we provide some basic notions and results about schemes. If the reader is well acquainted with them, then this section can be skipped and the reader can come back later when necessary. If the reader is not familiar with schemes, then this section should serve as a quick summary of the only basic notion and results that will be used in this paper. Reader more familiar with the classical language of varieties can also take a look at \cite{eisenbud2006geometry} or \cite[Appendix A.3-A.5]{greuel2008singular} for a comparison of classical and modern algebraic geometry. 

In this article, all rings are commutative with a multiplicative identity and all fields are of characteristic 0. 
\subsection{Affine Schemes}
Our main reference will be Hartshorne's classical textbook \emph{Algebraic Geometry} \cite{hartshorne2013algebraic}, supplemented by Liu's \emph{Algebraic Geometry and Arithmetic Curves} \cite{liu2002algebraic}. In case that a terminology is not explained here, please refer to \cite{hartshorne2013algebraic} or \cite{liu2002algebraic}. 



Let $S$ be a ring. The \textbf{spectrum} of $S$, denoted by $\spec S$, is a topological space together with a sheaf of regular functions on it. The topological space of $\spec S$ is the set of prime ideals in $S$, equipped with the topology that takes closed sets as $$V(I)=\{p\in\spec S|I\unlhd S,\ I\subseteq p\}.$$ This topology is named after Zariski. The \textbf{distinguished open sets} in $\spec S$ are the sets of the form
$$D(h)=\{p\in \spec S|h\in S,\ h\notin p\},$$ which is the complement of $V(\ideal{h})$. Distinguished open sets form a base for the Zariski topology \cite[Section II.2]{hartshorne2013algebraic}. Then $\mathbb{A}_S^n$, the \textbf{affine $n$-space} over $S$ is the spectrum of $S[x_1,\ldots,x_n]$, i.e.\ $\mathbb{A}_S^n=\spec S[x_1,\ldots,x_n]$. 

An \textbf{affine scheme} is the spectrum of some ring $S$. An affine $k$-scheme is the spectrum of a $k$-algebra. Therefore an \textbf{affine $k$-variety} is just the spectrum of a finitely generated $k$-algebra \cite[Definition 2.3.47]{liu2002algebraic}. Because we are primarily interested in affine schemes in this paper, and the category of affine schemes is equivalent to the category of commutative rings, with arrows reversed \cite[Prop. I-41]{eisenbud2006geometry}, we can often omit the structure sheaf.

A \textbf{morphism} $f:\spec B\to \spec A$ of affine schemes is given by the associated ring map $\varphi:A\to B$: suppose $p\in \spec B$, then $f(p)=\varphi^{-1}(p)$, which is a prime ideal in $A$.

A \textbf{fiber product} $X\times_Z Y$ of two affine schemes $X=\spec A$ and $Y=\spec B$ over a base scheme $Z=\spec S$ (i.e.\ there are morphisms $X\to Z$ and $Y\to Z$), is the spectrum of the tensor product $A\otimes_S B$.


If $X=\spec A$ is an affine scheme and $p\in X$, then the \textbf{residue field} of $p$ on $X$ is the residue field of the local ring at $p$: $\kappa_p= A_p/pA_p=\Frac (A/p)$ (see \cite[Exercise II.2.7]{hartshorne2013algebraic}). Suppose $X$ is an affine $k$-scheme, then $p\in X$ is said to be a \textbf{$k$-rational point} if the residue field of $p$ is $k$ (see \cite[Exercise II.2.8]{hartshorne2013algebraic}). The set of all the $k$-rational points in $X$ is denoted by $X(k)$.

Let $k$ be a field. We use the notation $k^n$ to denote the \textbf{classical affine $n$-space} (the $n$-fold cartesian product of $k$). Please be aware that $k^n\subsetneq \mathbb{A}_k^n$ and it inherits the Zariski topology from the affine $n$-spaces. In fact, if $X$ is an affine $k$-variety with a closed immersion $X\cong V(I)\subseteq \mathbb{A}_k^n$, then the set of the $k$-rational points $X(k)$ can be identified with the common solutions of $I$ in $k^n$ \cite[Proposition 3.2.18]{liu2002algebraic}. 
$$X(k) \xleftrightarrow{\text{ 1-1 }} \{(x_1,\ldots,x_n)\in k^n| f(x_1,\ldots,x_n)=0,\forall f\in I\}.$$
This bridges the gap between the classical algebraic geometry and the modern algebraic geometry. Now an algebraic set in $k^n$ is just the set of $k$-rational points of a $k$-variety. In particular, $k^n=\mathbb{A}_k^n(k)$.

Suppose $f:X\to Y$ is a $k$-morphism of affine $k$-schemes $X,Y$. We add the subscript $0$ to $f$ to denote the \textbf{restriction of $f$ on the $k$-rational points}: $$f_0:\begin{array}{rcl}X(k)&\to& Y(k)\\ x&\mapsto& f(x).\end{array}$$ If $k$ is algebraically closed and $X$, $Y$ are $k$-varieties, then $f_0$ is just the classical morphism of classical varieties.

If $f:X\to Y$ is a morphism of schemes, $y\in Y$, the \textbf{(scheme-theoretic) fiber} of $f$ over $y$ is defined to be the $\kappa_y$-scheme $$X_y=X\times_Y \spec \kappa_y.$$
Similarly, we define the \textbf{geometric fiber} of $f$ over $y$ to be 
$$X_{\overline{y}}=X_y\times_{\kappa_y} \spec\overline{\kappa_y}=X\times_Y \spec\kappa_y\times_{\kappa_y}\spec\overline{\kappa_y}= X\times_Y \spec\overline{\kappa_y},$$
where $\overline{\kappa_y}$ is the algebraic closure of $\kappa_y$.
One can show that the underlying topological space of $X_y$ is homeomorphic to the subset $f^{-1}(y)$ of $X$ (see \cite[Exercise II.3.10]{hartshorne2013algebraic}). When $X=\spec B$ and $Y=\spec A$ are affine, then $y$ is a prime ideal in $A$ so $X_y=\spec B\otimes_A \Frac(A/y)$ and $X_{\overline{y}}=\spec B\otimes_A \overline{\Frac(A/y)}$ by the definition of fiber product. 

A morphism $f:\spec B\to \spec A$ is a \textbf{finite} morphism, if $B$ is a finite $A$-module via the associated ring map. A morphism $f:X\to Y$ is \textbf{quasi-finite}, if $X_y=f^{-1}(y)$ is a finite set for every $y\in Y$. A finite morphism is always quasi-finite \cite[Exercise II.3.5]{hartshorne2013algebraic}.

Suppose $X=\spec B$ and $Y=\spec A$ are two affine schemes. Morphism $f:X\to Y$ is said to be \textbf{free}, if $B$ is a free $A$-module via the associated ring map $f^\#:A\to B$ (see \cite[Page 109]{hartshorne2013algebraic}). For readers familiar with flatness, a finite flat morphism is almost equivalent to a finite free morphism, because a finite module over a noetherian local ring is flat if and only if it is a free module \cite[Exercise 6.2]{eisenbud2013commutative}.

In a topological space $X$, a subset $S$ is said to be \textbf{locally closed}, if it is the intersection of an open subset and a closed subset. A \textbf{constructible set} is a finite union of locally closed sets \cite[Definition 10.7]{kemper2011course}. 

\subsection{Generic Freeness}

	Suppose $A$ is a Noetherian domain and $B$ is a finitely generated $A$-algebra. If $M$ is a finitely generated $B$-module, then there exists an element $0\ne a\in A$ such that $M_a$ is a free $A_a$-module. This is a classical result due to Grothendieck, called Grothendieck's Generic Freeness Lemma \cite[Lemme 6.9.2]{grothendieck1965elements}, which is more commonly known as ``Generic Flatness'' in algebraic geometry.
	
	The theory of Gr\"{o}bner basis gives it a constructive proof, which is known to many people as a folklore \cite{vasconcelos1997flatness}. We collect the results here.
\begin{lemma}[Effective Generic Freeness]\label{lem-specialization-of-GB}
	%
	Let $I\unlhd k[\bs{y},\bs{x}]$ and $J=I\cap k[\bs{y}]\unlhd k[\bs{y}]$, $T_x$ (resp. $T_y$) is a monomial order in $\bs{x}$ (resp. $\bs{y}$). $T=T_x>T_y$ is a block order. Suppose $\mathcal{G}$ is a reduced Gr\"{o}bner basis of $I$ with respect to $T$, $\mathcal{G}_y = \mathcal{G}\cap k[\bs{y}]$ and $\mathcal{G}_x=\mathcal{G}\backslash \mathcal{G}_y$. Then
	\begin{enumerate}[label*=(\alph*)]
		\item The set $\mathcal{G}_y$ is a Gr\"{o}bner basis of $J$ with respect to $T_y$; \label{item-elimination-property}
		\item For all $f\in \mathcal{G}_x$, the leading coefficient $\lc_{T_x}(f)\notin J$ \label{item-leading-coef-not-in-J};
		\item Set $w=\prod_{f\in \mathcal{G}_x}\lc_{T_x}(f)$, then for all prime ideals $p\unlhd k[\bs{y}]$ containing $J$: if $w\notin p$, then the specialization of $\mathcal{G}$ at $\kappa_p=\Frac(k[\bs{y}]/p)$ is a Gr\"{o}bner basis for $I\otimes_k \kappa_p$ with respect to $T_x$.
		\item Moreover, $k[\bs{y},\bs{x}]_w/I_w$ is a free $(k[\bs{y}]_w/J_w)$-module. A basis of this free module can be chosen to be the monomials in $\bs{x}$ not divisible by any leading monomial of elements of $\mathcal{G}_x$ with respect to $T_x$. \label{item-generic-freeness-lemma}
	\end{enumerate}
\end{lemma}
\begin{proof}~
	
	\begin{enumerate}[label*=(\alph*)]
		\item Elimination property of Gr\"{o}bner basis \cite[Section 3.1, Exercise 5]{cox2013ideals}.
		\item Since $\mathcal{G}$ is a reduced Gr\"{o}bner basis, no leading monomial can divide another leading monomial. If there is some $f\in \mathcal{G}_x$ such that $l=\lc_{T_x}(f)\in J$, then by the previous result, there is some $g\in \mathcal{G}_y$ such that $\lt_{T_y}(g)|\lt_{T_y}(l)$. Therefore $\lt_{T}(g) | \lt_T(f)$, which is a contradiction.
		\item Kalkbrener's stability criterion \cite[Theorem 3.1]{kalkbrener1997stability}.
		\item This is \cite[Theorem 10.1]{kemper2011course}.
	\end{enumerate}

\end{proof}

\section{Counting Geometry Fiber Cardinality}\label{sect-fib-cls-of-gen-finite-mor}
Our goal in this section is to count the cardinality of geometric fiber of a finite free morphism. Throughout this section, $k$ is a field of characteristic 0.

A function $\varphi: X\to \mathbb{Z}$ is said to be \textbf{lower} \textbf{semi-continuous}, if $\{x\in X|\varphi(x)\leq n\}$ is closed for all $n\in \mathbb{Z}$. The main result in this section is the lower semi-continuity of geometric fiber cardinality:

\begin{theorem}[Lower Semi-Continuity of Geometric Fiber Cardinality]\label{thm-gen-lsc-of-card-geo-fiber}
	
	If $f:X\to Y$ is a finite, free morphism between affine $k$-varieties $X$ and $Y$, then the cardinality of geometric fiber $y\mapsto \# X_{\overline{y}}$ is a lower semi-continuous function on $Y$.
	
	If $U=D(h)$ is a distinguished open subset of $X$, then $y\mapsto \# U_{\overline{y}}$ is also lower semi-continuous on $Y$.
\end{theorem}

 Now, we will introduce multivariate Hermite quadratic form, which is a trace pairing technique to count the number of geometric points of a zero-dimensional $k$-algebra. It will play a central role in our proof.
 The tool is proposed by Becker and W\"{o}rmann \cite{becker1991sums,becker1994trace}, and independently by Pedersen, Roy and Szpirglas \cite{pedersen1991counting,pedersen1993counting} in its current form. Scheja and Storch obtained essentially the same result in a different form \cite[Theorem 94.7]{scheja1988lehrbuch}. See \cite{cox2021stickelberger} for a discussion on the history. This technique is widely used in symbolic computation (for example, \cite{weispfenning1998new,perrucci2017elementary,le2022solving}). For our convenience, we will rephrase this in the language of schemes. 

We need a notation here. Suppose $S$ is a ring and $A$ is an $S$-algebra that is also a free $S$-module of finite rank. Let $f\in A$, then the \textbf{multiplication map} $L_f$ is the $S$-module endomorphism given by
$$L_f:\begin{array}{rcl}
	A &\to& A \\
	g&\mapsto& fg.
\end{array}$$
Because $A$ is a free $S$-module of finite rank, the trace of $L_f$ is well-defined as the sum of diagonal entries of a matrix representing $L_f$.

\begin{lemma}[Multivariate Hermite Quadratic Form]\label{lem-multi-HQF}
	Suppose $k$ is a field. Let $A$ be a zero-dimensional $k$-algebra. Let $\overline{A}=A\otimes_k \overline{k}$ and $X=\spec \overline{A}$. Define Hermite Quadratic Form on $A$:
	$$H: \begin{array}{rcl}A \times A &\to& k\\ (f,g)&\mapsto& \tr L_{fg},\end{array}$$
	then $\rank H=\# X$.
	
	Moreover, if $h\in A$, define
	$$H_h: \begin{array}{rcl}A \times A &\to& k\\ (f,g)&\mapsto& \tr L_{fhg},\end{array}$$
	then $\rank H_h=\# X_h$ ($X_h=\spec \overline{A_h}=\spec (A_h\otimes_k \overline{k})$).
	
	If $k$ is in addition ordered and the real closure of $k$ is $R$, let $A'=A\otimes_k R$ be the base change of $A$ to $R$ and $X'=\spec A'$, then $\sign H=\# X'(R)$, i.e., the signature of $H$ is exactly the number of $R$-rational points in $X'=\spec A'$.
	
	Similarly, $\sign H_h=\#\{p\in X'(R)|h(p)>0\}-\#\{p\in X'(R)|h(p)<0\}$.
\end{lemma}
\begin{proof}
	This is just \cite[Theorem 4.102]{basu2008algorithms}, rephrased in the language of schemes.
\end{proof}

\begin{corollary}[Relative Version]\label{cor-rel-multi-HQF}
	Let $k$ be a field of characteristic 0. Suppose $X=\spec B$, $Y=\spec A$ are two affine $k$-varieties, $p\in Y$ and $f:X\to Y$ is a $k$-morphism. Let $\kappa_p=\Frac A/p$ be the residue field of $p$. If the fiber $X_p$ or the geometric fiber $X_{\overline{p}}$ is finite, then the $\kappa_p$-bilinear form on $B\otimes_A \kappa_p$ sending $(f,g)$ to $\tr L_{fg}$ has its rank equal to $\#X_{\overline{p}}$. 
	
	Moreover, if $h\in B$ and $U=D(h)\subseteq \spec B$ is a distinguished open subset, then the $\kappa_p$-bilinear form on $B\otimes_A \kappa_p$ sending $(f,g)$ to $\tr L_{fhg}$ has its rank equal to $\#U_{\overline{p}}$. 
\end{corollary}

\begin{proof}
	Notice that the fiber $X_p$ is the spectrum of $B\otimes_A \kappa_p$ and the geometric fiber $X_{\overline{p}}$ is just the spectrum of $B\otimes_A \kappa_p \otimes_{\kappa_p} \overline{\kappa_p}$. Therefore everything follows from Lemma \ref{lem-multi-HQF}.
\end{proof}

%


In order to prove Theorem \ref{thm-gen-lsc-of-card-geo-fiber}, we extend the classical trace pairing technique to the free module of finite rank case now.
\begin{definition}\label{def-para-HQF}
	If $A$ is a finitely generated $k$-algebra and $B$ is an $A$-algebra that is also a free $A$-module of finite rank, then we define the \textbf{relative Hermite Quadratic Form} of $B$ to be:
	$$H:\begin{array}{rcl}
		B \times B & \to & A \\
		(f,g)  & \mapsto & \tr L_{fg}.
	\end{array}$$
	This is well-defined because $B$ is a free $A$-module of finite rank.
	
	Similarly, if $h\in B$, we define the \textbf{relative Hermite Quadratic Form associated with} $h$ of $B$ to be:
	$$H_h:\begin{array}{rcl}
		B \times B & \to & A \\
		(f,g)  & \mapsto & \tr L_{fhg}.
	\end{array}$$
	
\end{definition}
	
\begin{lemma}[Specialization Property of relative Hermite Quadratic Form]
	\label{lem-specialization-of-para-HQF}
	Notations as above and $p\in \spec A$, then the relative Hermite Quadratic Forms $H$ and $H_h$ evaluated at $p$ coincides with the classical Hermite Quadratic Form.
\end{lemma}
\begin{proof}
	Choose a basis $\{b_1,\ldots,b_r\}$ of $B$ as free $A$-module, then $H$ ($H_h$ resp.) can be represented as a matrix $\left(\tr L_{b_ib_j}\right)_{i,j}$ ($\left(\tr L_{b_ihb_j}\right)_{i,j}$ resp.). Let $p$ be a prime ideal of $A$ and $\varphi$ is the canonical map $\varphi: B\to B\otimes_A \Frac(A/p)$, the images of $b_i$ in $B\otimes_A A_p \otimes_{A_p} A_p/pA_p\cong B\otimes_A \Frac(A/p)$ is a basis for $\Frac(A/p)$-vector space $B \otimes_A \Frac(A/p)$.
	
	Now the matrix of the classical Hermite Quadratic Form is given by the traces of $L_{\varphi(b_ib_j)}$, which are exactly $\tr L_{b_ib_j}$ evaluated at $p$. The same argument works for $H_h$. The proof is completed.
\end{proof}

\begin{proof}[Proof of Theorem \ref{thm-gen-lsc-of-card-geo-fiber}]
	
	
	Let $X=\spec B$ and $Y=\spec A$, then $B$ is a free $A$-module of finite rank via the morphism, choose a basis $\{b_i\}$ for $B$.
	
	By Definition \ref{def-para-HQF} and Lemma \ref{lem-specialization-of-para-HQF}, the relative Hermite Quadratic Form specializes to all $p\in Y$, and Corollary \ref{cor-rel-multi-HQF} shows that the rank of Hermite Quadratic Form specialized at $p$ is equal to $\#X_{\overline{p}}$. Let $H=(\tr L_{b_ib_j})_{b_i,b_j\in \mathscr{B}}$, then $\rank H(p)< r$ if and only if all the $r$-minors of $H(p)$ vanish. Let $I_r$ be the ideal in $A$ generated by all $r$-minors of $H$, then $\{p\in Y|\rank H(p)\leq r\}=V(I_{r+1})$ is closed, i.e.\ $y\mapsto \# X_{\overline{y}}$ is lower semi-continuous on $Y$. By replacing $H$ with $H_h$, we can show the same for $y\mapsto \#U_{\overline{y}}$. This completes the proof.
\end{proof}

\begin{remark}
	Actually, Theorem \ref{thm-gen-lsc-of-card-geo-fiber} is a special case of a Grothendieck theorem \cite[Proposition 15.5.9]{grothendieck1966elements}, which states for a quasi-finite separated flat morphism locally of finite presentation between schemes, using non-constructive techniques like Stein factorization. 
\end{remark}

\section{The Geometry of Cylindrical Algebraic Decomposition}\label{sect-geometric-cad}
Throughout this section, $k$ is an ordered field and $R$ is the real closure of $k$. Let $C=R(i)=R[x]/\ideal{x^2+1}$ be the algebraic closure of $R$. 

We recall the following two basic facts in real algebraic geometry, which are about semi-algebraic connectedness and continuity of roots, respectively.

\begin{lemma}[{\cite[Proposition 3.12]{basu2008algorithms}}]\label{lem-locally-constant-function-is-constant}
	If $S$ is a semi-algebraically connected semi-algebraic set	and $f : S \to R$ is a locally constant semi-algebraic function (i.e.\ given $x \in S$,	there is an open $U\subset S$ such that for all $y \in U$, $f(y) = f(x)$), then $f$ is a constant.
\end{lemma}


\begin{theorem}[{\cite[Theorem 5.12]{basu2008algorithms}}, Continuity of Roots]\label{thm-continuity-of-roots}	
	Let $P\in R[x_1,\ldots,x_n]$ and let $S$ be a semi-algebraic subset of $R^{n-1}$. Assume that $\deg_{x_n}(P)$ is constant on $S$ and that for some $a'\in S$, $z_1,\ldots,z_j$ are the distinct roots of $P(a';x_n)$ in $C=R[i]$, with multiplicities $\mu_1,\ldots,\mu_j$, respectively. 
	
	If the open disks $B(z_i,r)\subset C$ are disjoint, then there is an open neighborhood $V$ of $a'$ such that for every $x'\in V\cap S$, the polynomial $P(x';x_n)$ has exactly $\mu_i$ roots, counted with multiplicities, in the disk $B(z_i,r)$, for $i=1,\ldots,j$.
\end{theorem}

\subsection{Geometric Delineability for a Variety}

In this subsection, we will extend Collins's notion of delineability to a variety. The goal is to show that a finite free morphism of $R$-varieties admits semi-algebraic sections on the regions that have a constant geometric fiber cardinality if it factors through the cylinder over the base. We will see that counting the geometric fiber lies at the heart of geometry of cylindrical algebraic decomposition.

\begin{lemma}[Rank-Signature Lemma]\label{lem-invariant-rank-implies-invariant-signature}
	If $S$ is a semi-algebraically connected semi-algebraic subset of $\mathop{\mathrm{Sym}}(R,n):=\{M\in R^{n\times n}| M_{ij}=M_{ji} \}$ (with the Euclidean topology inherited from $R^{n\times n}$) and all the matrices in $S$ are of the same rank, then all matrices in $S$ have the same signature.
\end{lemma}
\begin{proof}
	Define map 
	$$\chi:\begin{array}{rcl}
		S&\to& R^n \\
		M&\mapsto& (a_{n-1},a_{n-2},\ldots,a_0),
	\end{array}$$
	where $\lambda^n+\sum_{i=0}^{n-1} a_i \lambda^i$ is the characteristic polynomial of $M$. By some abuse of notations, we also use $\chi_M$ to denote the characteristic polynomial of $M$, since this would not cause any confusion. Because $M$ is symmetric and all the entries are in $R$, it is diagonalizable and all the eigenvalues are in $R$ \cite[Theorem 4.43]{basu2008algorithms}. Suppose all the matrices in $S$ are of the same rank $r$, then for all $M\in S$, there are exactly $r$ non-zero eigenvalues. 
	To show the signature is invariant, it suffices to show the number of positive eigenvalues (counted with multiplicity) is invariant now.
	
	 By Theorem \ref{thm-continuity-of-roots}, the continuity of roots in the coefficients, for all $M\in S$, there is a neighborhood $U\subseteq S$ of $M$ such that for all $M'\in U$, $\chi_M$ and $\chi_{M'}$ share the same number of positive roots, counted with multiplicity. Define function 
	 $$\psi:\begin{array}{rcl}
	 	S & \to & R \\
	 	M &\mapsto& \sum\limits_{\lambda>0} \nu_{\lambda}(\chi_M)
	 \end{array}$$
	that counts the number of positive roots with multiplicity, where $\nu_{\lambda}(\chi_M)$ is the order of vanishing of $\chi_M$ at $\lambda$. Then $\psi$ is a locally constant function. Because a locally constant semi-algebraic function defined on a semi-algebraically connected semi-algebraic set is a constant (see Lemma \ref{lem-locally-constant-function-is-constant}), it remains to show $\psi$ is a semi-algebraic function now.
	
	A polynomial $f$ has exactly $p$ distinct positive roots with multiplicity $d_1,\ldots,d_p$ is equivalent to 
	\begin{align*}
		(\exists x_1)\cdots(\exists x_p)\left(\left(\bigwedge_{i=1}^p x_i>0\right) \wedge \bigwedge_{i=1}^p\left(\bigwedge_{j=0}^{d_i-1}f^{(j)}(x_i)=0 \wedge f^{(d_i)}(x_i)\ne 0\right)\right.\\
		\wedge \left.(\forall y)\left(\left((y>0)\wedge (\bigwedge_{i=1}^p y\ne x_i)\right)\implies (f(y)\ne 0)\right)\right).
	\end{align*}
	So by the Tarski-Seidenberg principle \cite[Theorem 2.92]{basu2008algorithms}, 
	\begin{align*}
		C_{d_1,\ldots,d_p}:=\{(a_{n-1},\ldots,a_0)|\lambda^n+\sum_{i=0}^{n-1}a_i\lambda^i=0 \text{ has exactly }p\text{ positive roots}\\ \text{with multiplicity }d_1,\ldots,d_p \}
	\end{align*}
	is a semi-algebraic set. Then 
	\begin{align*}
		B_l:=\{(a_{n-1},\ldots,a_0)|\lambda^n+\sum_{i=0}^{n-1}a_i\lambda^i=0 \text{ has }l \text{ positive roots}\}\\
		=\bigcup_{d_1+\cdots+d_p=l}C_{d_1,\ldots,d_p}
	\end{align*} 
	is also semi-algebraic. So the graph of $\psi$ is just $\bigcup_{l=0}^n (\chi^{-1}(B_l)\times \{l\})$, which is clearly semi-algebraic. Therefore $\psi$ is a semi-algebraic function by definition. Hence $\psi$ is a constant function on $S$ and the signature of matrices in $S$ is invariant.

\end{proof}

By counting geometric fiber, we also learn information about the real fiber, as shown in the theorem below.

\begin{theorem}[Geometric Fiber v.s.\ Real Fiber]\label{thm-invariant-number-of-complex-solutions-implies-invariant-number-of-real-solutions}
	Suppose $X=\spec B, Y=\spec A$ are two $R$-varieties and $f:X\to Y$ is a finite free $R$-morphism. Let $h\in B$. Consider the induced morphism of $R$-rational points $f_0: X(R)\to Y(R)$. Let $S\subseteq Y(R)$ be a semi-algebraically connected semi-algebraic set such that $\# (D(h))_{\overline{y}}$ is a constant for all $y\in S$, then $\# (f_0^{-1}(y)\cap D(h))$ is a constant over $S$.
\end{theorem}
\begin{proof}
	Because we can define $H_{h^2}$, the relative Hermite Quadratic Form associated with $h^2$ of $B$ (Definition \ref{def-para-HQF}) and it specializes to all points in $\spec A$ (Lemma \ref{lem-specialization-of-para-HQF}). By Lemma \ref{lem-multi-HQF}, the rank of $H_{h^2}$ is invariant over $S$, since $D(h)=D(h^2)$. Then Lemma $\ref{lem-invariant-rank-implies-invariant-signature}$ tells us the signature of $H_{h^2}$ is also invariant over $S$. Therefore using Lemma \ref{lem-multi-HQF} again, we see that the cardinality of real fibers $\# (f_0^{-1}(y)\cap D(h))$ is invariant over $S$ (notice that $f_0^{-1}(y)=X_y(R)$).
\end{proof}


To show that the real fibers can be glued to sections, we need a stronger condition in morphism.
\begin{definition}
	Let $X$, $Y$ be affine $k$-varieties and $\pi: X\to Y$ is a $k$-morphism. Suppose $\pi$ factors through the cylinder $Y\times_k \mathbb{A}_k^1$:
	\begin{center}
		\begin{tikzcd}[column sep=small]
			X \arrow[rr, "\pi"]\arrow[dr,hook]& & Y \\
			& Y\times_k \mathbb{A}_k^1 \arrow[ur] & \\
		\end{tikzcd}
	\end{center}
	and $X\to Y\times_k \mathbb{A}_k^1$ is a closed immersion.
	Then we say $f:X\to Y$ is a \textbf{Cylindrical Projection}.
\end{definition}

The following lemma says a Cylindrical Projection is finite and free, if and only if $X$ can be cut out by a monic polynomial in the cylinder $Y\times_k \mathbb{A}_k^1=\mathbb{A}_Y^1$.
\begin{lemma}[{\cite[Prop. 4.1]{eisenbud2013commutative}}]\label{lem-f.g.-free-A-module-iff-principal}
	Let $A$ be a ring and let $J\lhd A[x]$ be an ideal in the polynomial ring in one variable over $A$. Let $B=A[x]/J$, and let $s$ be the image of $x$ in $B$.
	\begin{enumerate}[label=\alph*.]
		\item $B$ is generated by $\leq n$ elements as an $A$-module if and only if $J$ contains a monic polynomial of degree $\leq n$. In this case $B$ is generated by $1,s,\ldots,s^{n-1}$. In particular, $B$ is a finitely generated $A$-module if and only if $J$ contains a monic polynomial.
		\item $B$ is a finitely generated free $A$-module if and only if $J$ can be generated by a monic polynomial. In this case $B$ has a basis of the form $1,s,\ldots,s^{n-1}$.
	\end{enumerate}
\end{lemma}

The next theorem shows that, if $X\to Y$ is a finite free Cylindrical Projection of $R$-varieties, then the real fibers are continuous functions on the base.

\begin{theorem}\label{thm-freeness-implies-geometric-delineability}
	Suppose $\pi:X\to Y$ is finite free Cylindrical Projection between affine $R$-varieties. Let $S\subseteq Y(R)$ be a semi-algebraic connected semi-algebraic set on which $\#X_{\overline{y}}$ is a constant, then there are semi-algebraic continuous functions $\xi_1,\ldots,\xi_l:S\to \pi^{-1}_0(S)$ such that $\pi\circ \xi_i=\mathrm{id}_S$ and $\{\xi_i(y)\}_{i=1}^{l}$ is the real fiber of $y$ in $X(R)$.
\end{theorem}

\begin{proof}
	We may embed $Y$ in $\mathbb{A}_R^n=\spec R[y_1,\ldots,y_n]$ as a closed subscheme. Then by $\pi$ is a Cylindrical Projection, $X$ can be embedded in $\mathbb{A}_R^{n+1}=\spec R[y_1,\ldots,y_n,x]$ too.				
	Here is the diagram showing the morphisms among all the schemes we just mentioned:
	\begin{center}
		\begin{tikzcd}
			X  \arrow[d,"\pi"]\arrow[r,hook,"i_X","closed"'] & \mathbb{A}_{R}^{n+1}  \arrow[d,"\pi"]\\
			Y \arrow[r,hook,"i_Y","closed"'] & \mathbb{A}_R^n.
		\end{tikzcd}
	\end{center}
	And what we are proving is that there are sections $\xi_i$ making the induced diagram commute:
	\begin{center}
		\begin{tikzcd}
			\pi_0^{-1}(S) \arrow[r,phantom,"\subseteq"] \arrow[d,shift left=1px,"\pi_0"] & X(R) \arrow[r,hook,"{(i_X)}_0"] \arrow[d,"\pi_0"] & R^{n+1} \arrow[d,"\pi_0"] \\
			S \arrow[r,phantom,"\subseteq",near end]  \arrow[u,dashed, bend left=25,"\xi_i"] & Y(R) \arrow[r,"{(i_Y)}_0"] & R^n.
		\end{tikzcd}
	\end{center}

	Because $\pi$ is a Cylindrical Projection, we can choose a monic polynomial $f\in R[y_1,\ldots,y_n,x]$ such that $X$ is cut out by the locus of $f$ over the cylinder $Y\times_R \mathbb{A}^1_R$, using Lemma \ref{lem-f.g.-free-A-module-iff-principal}. By Theorem \ref{thm-invariant-number-of-complex-solutions-implies-invariant-number-of-real-solutions}, the number of real fibers is a constant $m$ on $S$. Define 
	$$\xi_i:\begin{array}{rcl}S&\to &\pi_0^{-1}(S)\\ \bs{y}&\mapsto& (\bs{y},x_i),\end{array}$$
	where $x_i$ is the $i$-th largest root of $f(\bs{y};x)=0$ ($i=1,\ldots,m$), then $\xi_i$ are semi-algebraic functions by definition.
	
	It remains to show the continuity. Choose some particular $\bs{y}\in S$ and let $x_1,\ldots,x_l$ be the distinct solutions of $f(\bs{y};x)=0$ in $C$ with multiplicities $\mu_1,\ldots,\mu_l$. By renumbering the roots if necessary, we may assume $x_1<\cdots<x_m$ are the real roots of $f(\bs{y};x)=0$. Suppose $U_1,\ldots,U_l$ are disjoint open Euclidean neighborhoods of $x_1,\ldots,x_l$ such that $U_{m+1}\cap R=\cdots=U_{l}\cap R=\varnothing$.
	

	By Theorem \ref{thm-continuity-of-roots}, the continuity of roots in the coefficients, there is a neighborhood $U$ in $S$ of $\bs{y}$ such that for all $\bs{y}'\in U$, the equation $f(\bs{y}';x)=0$ has exactly $\mu_i$ solutions in $U_i$ for all $1\leq i\leq l$ (counted with multiplicity). But $\# X_{\overline{y}}=l$ is invariant on $S$, forcing that there is exactly 1 solution of multiplicity $\mu_i$ to $f(\bs{y}';x)=0$ in $U_i$ for $1\leq i\leq l$. Also, for $m+1\leq i\leq l$, the root of $f(\bs{y}';x)=0$ in $U_i$ must be non-real because $U_i\cap R=\varnothing$. So for $1\leq i\leq m$, the root of $f(\bs{y}';x)=0$ in $U_i$ is real, hence equal to the last coordinate of $\xi_i(\bs{y}')$. This shows the continuity for $\xi_i$.

\end{proof}


To show our theory does generalize Collins's original cylindrical algebraic decomposition, we will show that \cite[Theorem 1]{collins1975quantifier} is a special case of Theorem \ref{thm-freeness-implies-geometric-delineability}.

\begin{definition}[Collins's Delineability]\label{def-collins-delineability}
	Let $f(x_1,\ldots,x_r)\in R[x_1,\ldots,x_r]\, (r\geq 2)$, and let $S\subseteq R^{r-1}$ be a semi-algebraic set. The roots of $f$ are delineable over $S$ and that $\xi_1,\ldots,\xi_k\, (k\geq 0)$ delineate the real roots of $f$ over $S$ if the following conditions are satisfied:
	\begin{enumerate}
		\item There are $m\geq k$ positive integers $\mu_i$ such that if $(a_1,\ldots,a_{r-1})\in S$ then $f(a_1,\ldots,a_{r-1},x)=0$ has exactly $m$ distinct roots (in $C$), with multiplicities $\mu_1,\ldots,\mu_m$.
		\item $\xi_1<\xi_2<\cdots<\xi_k$ are semi-algebraic continuous functions from $S$ to $R$.
		\item If $(a_1,\ldots,a_{r-1})\in S$ then $\xi_i(a_1,\ldots,a_{r-1})$ is a root of $f(a_1,\ldots,a_{r-1},x)=0$ of multiplicity $\mu_i$ for all $1\leq i\leq k$ and $(a_1,\ldots,a_{r-1})\in S$.
		\item If $(a_1,\ldots,a_{r-1})\in S$, $b\in R$ and $f(a_1,\ldots,a_{r-1},b)=0$ then $b=\xi_i(a_1,\ldots,a_{r-1})$ for some $i$.
	\end{enumerate}
\end{definition}
\begin{corollary}[Collins {cf.\cite[Theorem 1]{collins1975quantifier}}]\label{cor-collins-thm-1}
	Let $f(x_1,\ldots,x_r)$ be an $R$-polynomial, $r\geq 2$. Let $S$ be a semi-algebraically connected semi-algebraic subset of $R^{r-1}$. If $\lc_{x_r}(f)$ has no zero on $S$ and the number of distinct roots of $f$ in $x_r$ is invariant on $S$, then the roots of $f$ are delineable on $S$.
\end{corollary}
\begin{proof}
	Set $A=R[x_1,\ldots,x_{r-1}]$ and $B=R[x_1,\ldots,x_{r}]/\ideal{f}.$ Let $X=\spec B$ and $Y=\spec A$. There is a natural morphism $\pi:X\to Y$ induced by the ring map $R[x_1,\ldots,x_{r-1}]\to R[x_1,\ldots,x_{r}]/\ideal{f}$. If $\deg_{x_r}(f)=0$ then we are already done. So let us suppose $\deg_{x_r}(f)>0$ and set $w=\lc_{x_r}(f)$, then $B_w$ is a free $A_w$-module of finite rank (generated by monomials $1,x_r,\ldots,x_r^{\deg_{x_r}(f)-1}$). In other words, $\spec B_w\to \spec A_w$ is a finite free morphism. 
	The rest immediately follows from Theorem \ref{thm-freeness-implies-geometric-delineability} (note that the constant multiplicity is implied in the proof of Theorem \ref{thm-freeness-implies-geometric-delineability}.)
\end{proof}

Motivated by Collins's Delineability, we will define Geometric Delineability for a variety.
\begin{definition}[Geometric Delineability I]\label{def-geo-delineability-part-1}
	Let $T=V(f_1,\ldots,f_s)$ be a variety in $\mathbb{A}^{n+1}_R$, where $f_i\in R[x_1,\ldots,x_{n+1}]$. Let $S$ be a semi-algebraically connected semi-algebraic set in $R^n$. We say $T$ is \textbf{(strictly) geometrically delineable} over $S$, if the following conditions holds:
	\begin{enumerate}
		\item For all $a=(a_1,\ldots,a_n)\in S$, the number of distinct solutions (in $C$) of $$f_1(a;x_{n+1})=\cdots=f_s(a;x_{n+1})=0$$
		is finite and invariant.
		\item There are semi-algebraic continuous functions $\xi_1,\ldots,\xi_k: S\to T(R)$ such that $\pi\circ\xi_i=\mathrm{id}_S$, where $\pi: \mathbb{A}^{n+1}_R\to \mathbb{A}^{n}_R$ is the projection map. Additionally, we can require $\xi_i$ to be arranged in strictly ascending order in the last coordinate.
		\item If $t\in T(R)\cap (S\times R^1)$, then there is some unique $1\leq i\leq k$ such that $t=\xi_i(\pi(t))$.
	\end{enumerate}
	The semi-algebraic continuous functions $\xi_i$ are called the \textbf{sections} of $T$ over $S$.
	
	We make the convention that if $T\cap (S\times R^1)=S\times R^1$ (i.e., $T$ is a cylinder over $S$), then we say $T$ is \textbf{(vacuously) geometrically delineable} over $S$.
\end{definition}
\begin{example}[Geometric Delineability Examples]
	~
	\begin{figure}[hbtp]
		\centering
		\begin{subfigure}[t]{0.35\textwidth}
			\centering
			\includegraphics[width=\linewidth]{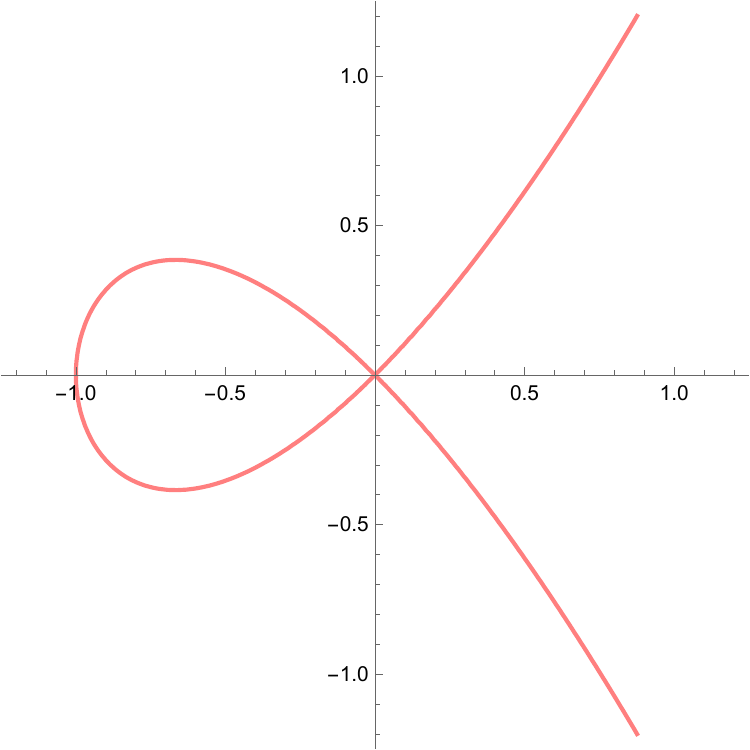}
			\caption{Nodal Curve}
			\label{fig:nodal-curve}
		\end{subfigure}
		~
		\begin{subfigure}[t]{0.35\textwidth}
			\centering
			\includegraphics[width=\linewidth]{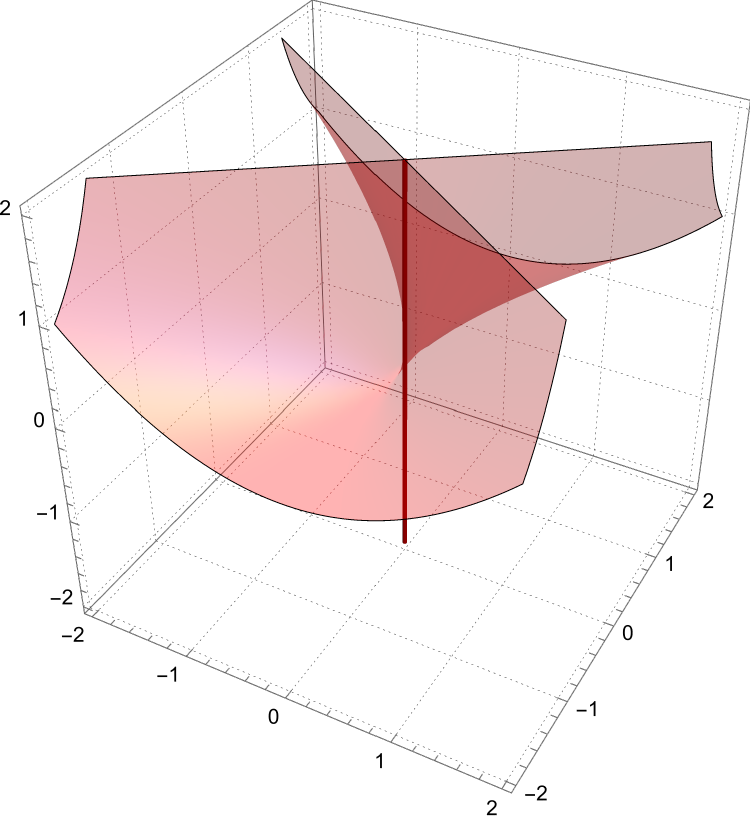}
			\caption{Whitney Umbrella}
			\label{fig:whitney-umbrella}
		\end{subfigure}	
		\caption{Several Geometric Objects to Demonstrate Geometric Delineability}
		\label{fig-example-2}
	\end{figure}

	Figure \ref{fig:nodal-curve} is the nodal curve $E: y^2-x^2(x+1)=0$ in $R^2$. Let $$I_1=(-\infty,-1), p_1=\{-1\}, I_2=(-1,0), p_2=\{0\} \text{ and } I_3=(0,+\infty)$$ be intervals or points in $R^1$. It is easy to see that $E$ is geometrically delineable over these subsets of $R^1$. For example, $E$ is geometrically delineable over $I_2$ because there are two different roots of $y^2-x^2(x+1)=0$ in $y$ when $x\in I_2$, given by $y=\pm\sqrt{x^2(x+1)}$.
	
	Figure \ref{fig:whitney-umbrella} is the Whitney's umbrella $U: x^2-y^2z=0$ in $R^3$. Clearly $U$ is geometrically delineable over $\{(x,y)\in R^2| y\ne 0\}$ because $z=\frac{x^2}{y^2}$ over this region. Also $U$ is geometrically delineable over $\{(x,y)\in R^2| x\ne 0,y=0\}$ because there is no solution in $z$ when $y=0$ but $x\ne 0$. By the convention we made in Definition \ref{def-geo-delineability-part-1}, $U$ is (vacuously) geometrically delineable over $\{(0,0)\}$ because the whole $z$-axis is contained in $U$.
	
\end{example}

We can rephrase Theorem \ref{thm-freeness-implies-geometric-delineability} in our new terminology by just saying $X$ is geometrically delineable over $S$.
\begingroup
\def\thetheorem{\ref{thm-freeness-implies-geometric-delineability}}
\begin{theorem}[Reformulation in Geometric Delineability]
	If $\pi:X\to Y$ is a finite free Cylindrical Projection between affine $R$-varieties. Let $S\subseteq Y(R)$ be a semi-algebraic connected semi-algebraic set such that $\# X_{\overline{y}}$ is invariant over $S$, then $X$ is geometrically delineable over $S$.
\end{theorem}
\addtocounter{theorem}{-1}
\endgroup
\begin{remark}
	While Collins's delineability in \cite{collins1975quantifier} requires constant multiplicity in each section, this requirement is dropped in our definition of geometric delineability. However, these two notions are equivalent for a finite free cylindrical projection from a hypersurface (recall the proof of Theorem \ref{thm-freeness-implies-geometric-delineability}). The constant multiplicity is never used in the c.a.d.\ construction, so we do not include it in the definition. There is an alternative definition of delineability from Arnon, Collins and McCallum that do not require a constant geometric fiber cardinality or a constant multiplicity \cite{arnon1984cylindrical}.
\end{remark}

We propose Algorithm \ref{algo-PROJ1} to compute the regions over which a variety $V(I)$ is geometrically delineable. The extra output will be useful later (Lemma \ref{lem-cad-set-operation}).

\begin{algorithm}[htbp]
	\KwIn{An ideal $I\unlhd k[\bs{y},x]$; Here, $k$ is a subfield of $R$.}
	\SetKwInput{KwOut}{Output 1}
	\KwOut{A list of pairs of ideals in $k[\bs{y}]$: $\{(J_i,J_i')\}$ such that $V(I)$ is geometrically delineable over any semi-algebraically connected semi-algebraic set contained in $V(J_i)\backslash V(J_i').$}
	\SetKwInput{KwOut}{Output 2}
	\KwOut{An extra list of ideals $\{\mathfrak{a}_i\}$ in $k[\bs{y}]$. For each tuples $(J_i,J_i')$ in the Output 1, the locally closed set $V(J_i)\backslash V(J_i')$ is a finite Boolean combination of the sets $\{V(\mathfrak{a}_i)\}$ (using $\cap,\cup,\backslash$).}
	\caption{UnramifiedStratification}\label{algo-PROJ1}
	
	Compute a reduced Gr\"{o}bner basis $\mathcal{G}$ of $I$ w.r.t some block order $x>T_{\bs{y}}$\;
	$\mathcal{G}_y:=\mathcal{G}\cap k[\bs{y}],\mathcal{G}_x:=\mathcal{G}\backslash\mathcal{G}_y $\;
	
	\If(\tcp*[f]{$V(I)$ is a cylinder}){$\mathcal{G}_x=\varnothing$}{\Return{$\{(\ideal{\mathcal{G}_y},\ideal{1})\},\{\ideal{\mathcal{G}_y}\}$\label{line-proj1-not-finite}} \tcp*{Vacuously geometrically delineable}}
	
	$w:=\prod_{f\in \mathcal{G}_x} \lc_{x}(f)$ \tcp*{Effective Generic Freeness}
	
	$H:=\mathtt{HermiteQuadraticForm}(\mathcal{G}), s:=\mathtt{sizeof}(H)$\; 
	
	$o_1:=\{(\ideal{\mathcal{G}_y}+\mathtt{minors}(i+1,H),\ideal{\mathcal{G}_y}+\mathtt{minors}(i,H)\cdot \ideal{w})\}_{i=0}^s$
	\tcc*{$\mathtt{minors}(i,H)$ computes the ideal generated by all $i$-minors in $H$. It returns $\ideal{1}$ when $i=0$, and $\ideal{0}$ when $i=s+1$.}
	$o_2:=\{\ideal{\mathcal{G}_y}+\mathtt{minors}(i+1,H)\}_{i=0}^s$\;
	
	\ForEach(\tcp*[f]{Recursively look at the non-free locus}){$f\in \mathcal{G}_x$}{
		$w_0:=\lc_{x}(f)$\;
		\If{$\ideal{\mathcal{G}_y}+\ideal{w_0}\neq (I+\ideal{w_0})\cap k[\bs{y}]$}{
			$o_2:=o_2\cup \{\ideal{\mathcal{G}_y}+\ideal{w_0}\}$\; \label{line-proj1-not-dorminant}
		}
		Merge $o_1$ and $o_2$ with $\mathtt{UnramifiedStratification}(I+\ideal{w_0})$\;		
	}
	
	\Return{$o_1$, $o_2$}\;
	
\end{algorithm}

\begin{algorithm}[htbp]
	\KwIn{A reduced Gr\"{o}bner Basis $\mathcal{G}$ of an ideal $I\unlhd k[\bs{y},x]$ with respect to some block order $x>T_{\bs{y}}$.}
	\KwOut{The relative Hermite Quadratic Form $H$.}
	\caption{HermiteQuadraticForm}\label{algo-Hermite-Quadratic-Form}
	
	$\mathcal{G}_y:=\mathcal{G}\cap k[\bs{y}],\mathcal{G}_x:=\mathcal{G}\backslash\mathcal{G}_y $\;
	Choose $f\in \mathcal{G}_x$ of the minimal degree, set $w=\mathtt{LC}_x(f)$\;
	$s:=\deg f$ \tcp*{The monomial basis $\mathscr{B}=\{1,x,\ldots,x^{s-1}\}$}
	\For{$i:=0$ \KwTo $s-1$}{
		\For(\tcp*[f]{Compute the trace of $L_{x^{i+j}}$}){$j:=0$ \KwTo $s-1$}{ 
			$H_{i,j}:=0$\;
			\For(\tcp*[f]{Compute the image of $x^k$ under $L_{x^{i+j}}$}){$k:=0$ \KwTo $s-1$}{
				$r:=\mathtt{PolynomialDivisionRemainder}(x^{i+j+k},f)$
				\tcp*{Polynomial Division over $k(\bs{y})[x]$. This does the NormalForm computation}
				$H_{i,j}:=H_{i,j}+\mathtt{coef}(r,k)$ \tcp*{the $k$-th coefficient of $r$}
			}
		}
	}
	Multiply $H$ by a suitable power of $w$ so entries of $H$ are in $k[\bs{y}]$;\\
	\tcc{The entries of $H$ are in $k[y]_w$, so there is some power $w^\ell$ that cancels all the denominators. This does not affect the rank of $H$.}
	\Return{$H$}\;
\end{algorithm}

\begin{remark}
	The reason for this name ``\texttt{UnramifiedStratification}'' is that the cardinality of geometric fiber remains the same over the regions defined by each pair of ideals in Output 1, so the fiber does not ``ramify'' into different branches.
	It should be also noted that, the regions defined by Output 1 are not necessarily disjoint.
\end{remark}

\begin{theorem} \label{thm-correctness-of-algo-PROJ1}
	Algorithm \ref{algo-PROJ1} terminates. For each tuples $(J_i,J_i')$ in the Output 1, $V(I)$ is geometrically delineable over any semi-algebraic connected semi-algebraic set contained in $V(J_i)\backslash V(J_i')$.
	
	Moreover, for each tuples $(J_i,J_i')$ in the Output 1, the locally closed set $V(J_i)\backslash V(J_i')$ can be obtained by finite intersections, unions and differences of varieties defined by the ideals $\{\mathfrak{a}_i\}$ in the Output 2. 
\end{theorem}
\begin{proof}
	\textit{Termination}. It suffices to show that there is no infinite recursion. By Lemma \ref{lem-specialization-of-GB}\ref{item-leading-coef-not-in-J}, for every $f\in \mathcal{G}_x$, the leading coefficient $w_0$ is not in $I$. So $I+\ideal{w_0}$ is a strictly larger ideal. Since $k[\bs{y},x]$ is Noetherian, the algorithm certainly terminates.
%

	\textit{Correctness}. The correctness can be proved by a standard technique called ``Noetherian Induction'' \cite[Exercise II.3.16]{hartshorne2013algebraic}. In other words, it suffices to show that if our algorithm works for closed subschemes of $V(I)$, then it works for $V(I)$ too. 
	
	In Line \ref{line-proj1-not-finite}, the algorithm finds that $I$ is fully contained in $k[\bs{y}]$, so $V(I)$ is vacuously geometrically delineable over $V(\ideal{\mathcal{G}_y})$ and our algorithm is correct in this case.
	
	Now suppose $\mathcal{G}_x\neq \varnothing$, then by Effective Generic Freeness (Lemma \ref{lem-specialization-of-GB}\ref{item-generic-freeness-lemma}), $k[\bs{y},x]_w/I_w$ is a free $(k[\bs{y}]_w/J_w)$-module, where $J=\ideal{\mathcal{G}_y}=I\cap k[\bs{y}]$. Also, this module is necessary of finite rank, otherwise $\mathcal{G}_x$ cannot be empty. Therefore define $X=(\spec k[\bs{y},x]_w/I_w)\times_k R$ and $Y=(\spec k[\bs{y}]_w/J_w)\times_k R$, the morphism $$\pi: X\to Y $$ is a finite free Cylindrical Projection of $R$-varieties.
	
	Next, the algorithm computes relative Hermite Quadratic Form $H$ by calling Algorithm \ref{algo-Hermite-Quadratic-Form}. The correctness of Algorithm \ref{algo-Hermite-Quadratic-Form} is quite direct, as it simply computes the matrix of left multiplication by $x^{i+j}$ using the monomial basis $\mathscr{B}$ and take the trace.
	
	Then, the algorithm calculates ideals of minors. By the specialization property of relative Hermite Quadratic Form (Lemma \ref{lem-specialization-of-para-HQF}) and Corollary \ref{cor-rel-multi-HQF}, if all $(i+1)$-minors vanish at $p\in Y$, but some $i$-minor does not vanish at $p$, then the rank of $H$ at $p$ is $i$ and $\#X_{\overline{p}}=i$.
	
	As a consequence, the geometric fiber cardinality is a constant $i$ on 
	\begin{equation}\label{eqn-geometrically-delineable-regions}
		\begin{array}{rl}
			&D(w)\cap \left(V\left(\ideal{\mathcal{G}_y}+\mathtt{minors}(i+1,H)\right)\backslash V\left(\ideal{\mathcal{G}_y}+\mathtt{minors}(i,H)\right)\right)\allowbreak\\
			=&V\left(\ideal{\mathcal{G}_y}+\mathtt{minors}(i+1,H)\right)\backslash V\left(\ideal{\mathcal{G}_y}+\mathtt{minors}(i,H)\cdot \ideal{w}\right)\\
			=&V\left(\ideal{\mathcal{G}_y}+\mathtt{minors}(i+1,H)\right)\backslash \left(V\left(\ideal{\mathcal{G}_y}+\mathtt{minors}(i,H)\right)\cup V(\ideal{\mathcal{G}_y}+\ideal{w})\right).
		\end{array}
	\end{equation} 
	Suppose $S$ is a semi-algebraically connected semi-algebraic set contained in it, then by Theorem \ref{thm-freeness-implies-geometric-delineability}, we can conclude that $V(I)$ is geometrically delineable over $S$. And the set in Equation (\ref{eqn-geometrically-delineable-regions}) can be obtained from finite intersections, unions and differences of the varieties defined by Output 2, if $V(\ideal{\mathcal{G}_y}+\ideal{w})$ is so. 
	
	Finally, Algorithm \ref{algo-PROJ1} takes care of the non-free locus $V(\ideal{\mathcal{G}_y}+\ideal{w})$, which is the union $\bigcup_{f\in \mathcal{G}_x} V(\ideal{\mathcal{G}_y}+\ideal{\lc_{x}(f)}) $. For any $f\in\mathcal{G}_x$, let $w_0=\lc_{x}(f)$, $I'=I+\ideal{w_0}$ and $J'=I'\cap k[\bs{y}]$. Clearly, if $V(I')$ is geometrically delineable over some $S\subseteq V(J')$ then $V(I)$ is also geometrically delineable over $S$ because $V(I')=V(I)\cap (V(w_0)\times \mathbb{A}_R^1)$. So the recursion is correct. Also, Line \ref{line-proj1-not-dorminant} guarantees that $V(\ideal{\mathcal{G}_y}+\ideal{w})$ is a union of varieties defined by Output 2.

\end{proof}

\subsection{Geometric Delineability for Several Varieties}

We now extend Geometric Delineability to multiple varieties, just as Collins extended his delineability to multiple polynomials.
\begin{definition}[Geometric Delineability II]\label{def-geo-delineability-part-2}
	A family of varieties $\{V_i\}_{i=1}^l$ in $\mathbb{A}^{n+1}_R$ is said to be \textbf{(jointly) geometrically delineable} over a semi-algebraically connected semi-algebraic set $S\subseteq R^n$, if each $V_i$ is geometrically delineable over $S$, and all the sections of $V_i$ ($i=1,\ldots,l$) are either identical or disjoint. As a consequence, the set of sections $\{\xi_{ij}:S\to R^{n+1}\}$ can be sorted in strictly ascending order in the last coordinate.
\end{definition}

The next lemma concludes the joint geometric delineability from the geometric delineability of the intersection.

\begin{lemma}\label{lem-joint-geo-del-from-intersection}
	Let $X$, $Y$ be two varieties in $\mathbb{A}_R^{n+1}$, and $Z=X\cap Y$. Suppose $X$, $Y$ and $Z$ are all geometrically delineable over a semi-algebraically connected semi-algebraic set $S\subseteq R^n$. Then $\{X,Y\}$ is jointly geometrically delineable over $S$.
\end{lemma}
\begin{proof}
	If either $X$ or $Y$ is vacuously geometrically delineable over $S$ then the lemma is trivially true and we can simply rule out this case. 
	
	Now suppose $\xi$ ($\eta$ resp.) is a section of $X$ ($Y$ resp.) over $S$ and $\xi(s_0)=\eta(s_0)$ for some $s_0\in S$.
	Since $Z=X\cap Y$ is geometrically delineable over $S$, there is a section $\zeta$ of $Z$ such that $\xi(s_0)=\eta(s_0)=\zeta(s_0)$. Because $Z\subseteq X$, for every $s\in S$, $\zeta(s)$ has to agree with some section of $X$. Then by the continuity, $\zeta$ has to agree with the same section in an Euclidean neighborhood of $s$. Therefore $$S_\xi=\{s\in S| \zeta(s)=\xi(s)\}$$ is an Euclidean open non-empty semi-algebraic subset of $S$. By an identical argument, we can show that for every section $\xi^*$ of $X$: $\{s\in S| \zeta(s)=\xi^*(s)\}$ is an Euclidean open semi-algebraic subset of $S$. So these open subsets constitute a finite disjoint open cover for $S$. But $S$ is semi-algebraically connected, forcing all but one are them is empty. This shows $S_\xi=S$, i.e.\ $\zeta$ and $\xi$ agree everywhere. Similarly $\zeta$ and $\eta$ must agree on $S$. Therefore all sections of $X$ and $Y$ are either identical or disjoint. 
\end{proof}

Therefore, to show the joint geometric delineability of a family of varieties $\{V_i\}$, it suffices to show the geometric delineability of each variety $V_i$ and intersection $V_i\cap V_j$ for each pair of varieties. The next theorem shows that we can conclude the geometric delineability of a closed subscheme $Z\subseteq X$ (an intersection for example) from a weaker condition, if we have the geometric delineability of $X$.

\begin{theorem}\label{thm-geo-del-of-closed-subscheme}
	Suppose $\pi:X\to Y$ is a finite free Cylindrical Projection of $R$-varieties. Let $Z$ be a closed subscheme of $X$. If $\pi|_Z:Z\to Y$ is also finite free and $S\subseteq Y(R)$ is a semi-algebraically connected semi-algebraic set satisfying that $\# X_{\overline{y}}$ is a constant for all $y\in S$, then $\# Z_{\overline{y}}$ is also a constant for all $y\in S$.
	
	As a corollary, $Z$ is geometrically delineable over $S$.
\end{theorem}
\begin{proof}
	Let $Y=\spec A$. Since $\pi$ is a finite free morphism, $X$ can be identified as $\spec A[x]/\ideal{f}$ for some monic $f\in A[x]$ by Lemma \ref{lem-f.g.-free-A-module-iff-principal}. Because $\pi|_Z: Z\to Y$ is finite and free, $Z\cong \spec A[x]/\ideal{g}$ for some monic $g\in A[x]$ by Lemma \ref{lem-f.g.-free-A-module-iff-principal} again.
	
	Now embed $Y$ in some affine space $\mathbb{A}_R^n=\spec R[y_1,\ldots,y_n]$ as a closed subscheme, so the monic polynomials $f,g$ can be taken from $R[y_1,\ldots,y_n,x]$ by a slight abuse of notations. Let $(c_1,\ldots,c_n)\in S$ and let $x_1,\ldots,x_l\in C$ be the distinct solutions of $f(c_1,\ldots,c_n;x)=0$ in $C$. Denote the multiplicity of $x_i$ as a root of $f(c_1,\ldots,c_n;x)=0$ by $\mu_i$. Notice that $Z$ is a subscheme of $X$, so by renumbering the roots if necessary, we can assume that $x_1,\ldots,x_m$ are the all distinct solutions of $g(c_1,\ldots,c_n;x)=0$ in $C$ ($0\leq m\leq l$). Denote the multiplicity of $x_i$ as a root of $g(c_1,\ldots,c_n;x)=0$ by $\lambda_i$.
	
	Choose disjoint open Euclidean neighborhoods $U_1,\ldots,U_l$ for $x_1,\ldots,x_l$. By the continuity of roots (Theorem \ref{thm-continuity-of-roots}), there is a open Euclidean neighborhood $U$ of $(c_1,\ldots,c_n)$ such that for all $(c'_1,\ldots,c'_n)\in U\cap S$, the equation $f(c'_1,\ldots,c'_n;x)=0$ has exactly $\mu_i$ solutions, counting multiplicity in $U_i$ for all $1\leq i\leq l$, and the equation $g(c'_1,\ldots,c'_n;x)=0$ has exactly $\lambda_i$ solutions, counting multiplicity in $U_i$ for all $1\leq i\leq m$.
	
	Because $\# X_{\overline{y}}$ is a constant over $S$, for all $(c'_1,\ldots,c'_n)\in U\cap S$, the equation $f(c'_1,\ldots,c'_n;x)=0$ should have only one solution of multiplicity $\mu_i$ in $U_i$ for all $1\leq i \leq l$. Notice that the roots of $g(c'_1,\ldots,c'_n;x)=0$ are also roots of $f(c'_1,\ldots,c'_n;x)=0$, hence the equation $g(c'_1,\ldots,c'_n;x)=0$ have only one solution of multiplicity $\lambda_i$ in $U_i$ for all $1\leq i\leq m$, and these are all the roots of $g$ by counting degree. Therefore for all $(c'_1,\ldots,c'_n)\in U\cap S$, the number of distinct complex solutions to the equation $g(c'_1,\ldots,c'_n;x)=0$ is a constant $m$. In other words, $\# Z_{\overline{y}}$ is a locally constant function on $S$. It is easy to show that this is a semi-algebraic function using an argument similar to Lemma \ref{lem-invariant-rank-implies-invariant-signature}. So by Lemma \ref{lem-locally-constant-function-is-constant}, $\# Z_{\overline{y}}$ is a constant on $S$.
	
	The last assertion follows immediately from Theorem \ref{thm-freeness-implies-geometric-delineability}.
\end{proof}

\begin{corollary}\label{cor-geo-del-of-closed-subscheme}
	Suppose $X,Z$ are closed subschemes of $\mathbb{A}_R^{n+1}$ and $Z$ is a closed subscheme of $X$. Let $Y,W$ be affine subschemes of $\mathbb{A}_R^n$ and let $\pi:\mathbb{A}_R^{n+1}\to \mathbb{A}_R^n$ be the projection map (a subscheme is an open subscheme of a closed subscheme). If $X\times_{\mathbb{A}_R^n} Y\to Y$ and $Z\times_{\mathbb{A}_R^n} W\to W$ are finite free morphisms, and $S\subseteq (Y\cap W)(R)$ is a semi-algebraically connected semi-algebraic set on which $X_{\overline{y}}$ is a constant, then $Z$ is geometrically delineable over $S$.
\end{corollary}
\begin{proof}
	We have the following commutative diagram.
	\begin{center} 
		\begin{tikzcd}[scale cd=1.0]
				& Z\times_{\mathbb{A}_R^n} (Y\cap W) \arrow[dl,hook,"closed"{rotate=18,above}]\arrow[dr,phantom,"\square"]\arrow[r] \arrow[d] & Z\times_{\mathbb{A}_R^n} W \arrow[d] \\
			X\times_{\mathbb{A}_R^n} (Y\cap W) \arrow[r]\arrow[dr,phantom,"\square"]\arrow[d] & Y\cap W \arrow[r]\arrow[dr,phantom,"\square"]\arrow[d] & W \arrow[d] \\
			X\times_{\mathbb{A}_R^n} Y \arrow[r] & Y \arrow[r] & \mathbb{A}_R^n
		\end{tikzcd}
	\end{center}
	Notice that finite free morphisms are stable under base change, so $Z\times_{\mathbb{A}_R^n} (Y\cap W)\to Y\cap W$ and $X\times_{\mathbb{A}_R^n} (Y\cap W)\to Y\cap W$ are also finite free morphisms. Applying Theorem \ref{thm-geo-del-of-closed-subscheme} to the upper-left triangle in the diagram completes the proof.
\end{proof}


We present Algorithm \ref{algo-PROJ2}, which stratifies the image of a projection into free-loci.
\begin{algorithm}[htbp]
\KwIn{An ideal $I\unlhd k[\bs{y},x]$; Here, $k$ is a subfield of $R$.}
\SetKwInput{KwOut}{Output 1}
\KwOut{A list of pairs of ideals in $k[\bs{y}]$: $\{(J_i,J_i')\}$ such that the projection is free on $V(J_i)\backslash V(J_i')$}
\SetKwInput{KwOut}{Output 2}
\KwOut{An extra list of ideals $\{\mathfrak{a}_i\}$ in $k[\bs{y}]$. For each tuples $(J_i,J_i')$ in the Output 1, the locally closed set $V(J_i)\backslash V(J_i')$ is a finite Boolean combination of the sets $\{V(\mathfrak{a}_i)\}$ (using $\cap,\cup,\backslash$)}
\caption{FreeStratification}\label{algo-PROJ2}

Compute a reduced Gr\"{o}bner basis $\mathcal{G}$ of $I$ w.r.t some block order $x>T_{\bs{y}}$\;
$\mathcal{G}_y:=\mathcal{G}\cap k[\bs{y}],\mathcal{G}_x:=\mathcal{G}\backslash\mathcal{G}_y $\;

\If(\tcp*[f]{$V(I)$ is a cylinder}){$\mathcal{G}_x=\varnothing$}{\Return{$\{(\ideal{\mathcal{G}_y},\ideal{1})\},\{\ideal{\mathcal{G}_y}\}$\label{line-proj2-not-finite}} \tcp*{Clearly a free morphism}}

$w:=\prod_{f\in \mathcal{G}_x} \lc_{x}(f)$ \tcp*{Effective Generic Freeness}


$o_1:=\{(\ideal{\mathcal{G}_y},\ideal{\mathcal{G}_y}+ \ideal{w})\}$;\\
\tcc{Notice that $V(\ideal{\mathcal{G}_y})\backslash V(\ideal{\mathcal{G}_y}+ \ideal{w})=\spec k[\bs{y}]_w/\ideal{\mathcal{G}_y}_w$ is affine. Thus it makes sense to say the projection is free on this set.}
$o_2:=\{\ideal{\mathcal{G}_y}\}$\;

\ForEach(\tcp*[f]{Recursively look at the non-free locus}){$f\in \mathcal{G}_x$}{
	$w_0:=\lc_{x}(f)$\;
	\If{$\ideal{\mathcal{G}_y}+\ideal{w_0}\neq (I+\ideal{w_0})\cap k[\bs{y}]$}{
		$o_2:=o_2\cup \{\ideal{\mathcal{G}_y}+\ideal{w_0}\}$\; \label{line-proj2-not-dorminant}
	}
	Merge $o_1$ and $o_2$ with $\mathtt{FreeStratification}(I+\ideal{w_0})$\;		
}

\Return{$o_1$, $o_2$}\;
\end{algorithm}

\begin{theorem} \label{thm-correctness-of-Proj2}
	Algorithm \ref{algo-PROJ2} terminates. For each tuples $(J_i,J_i')$ in the Output 1, the projection is free over $V(J_i)\backslash V(J_i')$.
	
	Moreover, for each tuples $(J_i,J_i')$ in the Output 1, the locally closed set $V(J_i)\backslash V(J_i')$ can be obtained by finite intersections, unions and differences of varieties defined by the ideals $\{\mathfrak{a}_i\}$ in the Output 2. 
\end{theorem}
\begin{proof}
	The proof is identical to Theorem \ref{thm-correctness-of-algo-PROJ1}. 	
\end{proof}

\subsection{The Algorithm}
We assume that field operations in $R$ ($+,-,\times,/,<,=$) can be effectively performed, so the real roots of a univariate equation can be effectively represented and computed (e.g.\ $R=\mathbb{R}_{\mathrm{alg}}$). Solving equations can be done by the real root isolation algorithm when $R$ is archimedean (see \cite[Algorithm 10.55]{basu2008algorithms}) and Thom Encoding if $R$ is a general real closed field (see \cite[Algorithm 10.105, 10.106]{basu2008algorithms}).


The following simple lemma is crucial in our algorithm for cylindrical algebraic decomposition.

\begin{lemma}\label{lem-cad-set-operation}
	Let $\mathcal{T}=\{T_i\}$ be a finite family of semi-algebraic sets in $R^n$, and $\mathfrak{D}$ is a cylindrical algebraic decomposition adapted to $\mathcal{T}$. Suppose $\mathcal{T}'=\{T'_j\}$ is another finite family of semi-algebraic sets in $R^n$ such that each $T'_j$ can be obtained by a finite number of intersections, unions, and differences of elements in $\mathcal{T}$. Then $\mathfrak{D}$ is also a cylindrical algebraic decomposition adapted to $\mathcal{T}'$.
\end{lemma}
\begin{proof}
	It suffices to prove that for any $T_i,T_j\in \mathcal{T}$, $\mathfrak{D}$ is also adapted to $\mathcal{T}\cup \{T_i\cap T_j, T_i\cup T_j, T_i\backslash T_j\}$.
	
	To see this, write $T_i=\cup_{\alpha\in \Lambda} S_{\alpha}$ and $T_j=\cup_{\alpha\in \Lambda'} S_{\alpha}$ as unions of disjoint cells of level $n$. Then $T_i\cup T_j=\cup_{\alpha\in \Lambda\cup \Lambda'} S_{\alpha}$ is a union of cells of level $n$. Similarly $T_i\cap T_j=\cup_{\alpha\in \Lambda\cap \Lambda'} S_{\alpha}$ and $T_i\backslash T_j=\cup_{\alpha\in \Lambda\backslash \Lambda'} S_{\alpha}$.
\end{proof}

\begin{algorithm}[htbp]
	\KwIn{A list of ideals $\{I_k\}_{k=1}^s$ in $k[\bs{y},x]$.}
	\SetKwInput{KwOut}{Output 1}
	\KwOut{A list of pairs of ideals in $k[\bs{y}]$: $\{(J_i,J_i')\}$.}
	\SetKwInput{KwOut}{Output 2}
	\KwOut{An extra list of ideals $\{\mathfrak{a}_i\}$ in $k[\bs{y}]$.}	
	\caption{CADProjection}\label{algo-CAD-Projection}
	$o_1=\varnothing, o_2=\varnothing$\;
	\For{$k:=1$ \KwTo $s$}{
			Merge $o_1,o_2$ with $\mathtt{UnramifiedStratification}(I_k)$
		}
	\For{$k:=1$ \KwTo $s$}{
			\For{$j:=1$ \KwTo $k-1$}{Merge $o_1,o_2$ with $\mathtt{FreeStratification}(I_j+I_k)$\;}
		}
	\Return{$o_1,o_2$}\;
\end{algorithm}

Algorithm \ref{algo-CAD-Projection} reduces a higher-dimensional c.a.d.\ computation to a lower-dimensional c.a.d.\ computation. It can be thought as an adaptation of Projection Operator to the geometric context.

\begin{theorem}\label{thm-CAD-lifting}
	Suppose $\{V(I_k)\}$ is a finite family of varieties in $\mathbb{A}_R^n$. A cylindrical algebraic decomposition $\mathfrak{D}'$ adapted to varieties defined by Output 2 of Algorithm \ref{algo-CAD-Projection}, $\mathtt{CADProjection}(\{I_k\})$ can be lifted to a cylindrical algebraic decomposition $\mathfrak{D}$ adapted to $\{V(I_k)\}$.
\end{theorem}
\begin{proof}
	Observe that cells in a cylindrical algebraic decomposition are always semi-algebraically connected, because a cell of level $>1$ is either a graph of a semi-algebraic continuous function defined over its silhouette, or a band between two graphs and a cell of level $1$ is either an interval or a point, so our assertion follows from an easy induction. 
	
	Using Lemma \ref{lem-cad-set-operation}, it suffices to prove that a cylindrical algebraic decomposition $\mathfrak{D}'$ adapted to the family of constructible sets $\{V(J_i)\backslash V(J_i')\}$ can be lifted to a cylindrical algebraic decomposition $\mathfrak{D}$ of $\{V(I_k)\}$, where $\{(J_i,J_i')\}$ is the Output 1 (recall Theorem \ref{thm-correctness-of-algo-PROJ1} and Theorem \ref{thm-correctness-of-Proj2}).
	
	Let $S$ be a level $n-1$ cell of $\mathfrak{D}'$. For each $I_k$, since $\mathtt{CADProjection}$ contains the Output 1 of $\mathtt{UnramifiedStratification}(I_k)$, $V(I_k)$ is geometrically delineable over $S$ by Theorem \ref{thm-correctness-of-algo-PROJ1}.
	
	Next, we will show the joint geometric delineability of every pair of varieties $(V(I_j),V(I_{k}))$ on $S$. By Lemma \ref{lem-joint-geo-del-from-intersection}, it suffices to show the intersection $V(I_j)\cap V(I_{k})=V(I_j+I_{k})$ is geometrically delineable over $S$. We may assume that $V(I_j)$ is not a cylinder over $S$ and $V(I_j+I_k)$ is non-empty over $S$ (otherwise the geometric delineability is trivial). Because $\mathtt{CADProjection}$ contains the Output 1 of $\mathtt{FreeStratification}(I_j+I_k)$, the projection of $V(I_j+I_k)$ is free over $S$. By applying Corollary \ref{cor-geo-del-of-closed-subscheme}, we see that $V(I_j+I_k)$ is indeed geometrically delineable over $S$.
	
	It remains to write down a cylindrical algebraic decomposition $\mathfrak{D}$ adapted to $\{V(I_k)\}$. We start with $\mathfrak{D}'$ and construct the cells of level $n$ to complete $\mathfrak{D}$. Let $S$ be cell of level $n-1$ in $\mathfrak{D}'$, each $V(I_k)$ is geometrically delineable over $S$ and the sections are either disjoint or identical. So we let the image of sections to be the cells of level $n$ over $S$. If there is some $V(I_k)$ vacuously geometrically delineable over $S$ then we also let the bands between sections/above the highest section/below the lowest section to be the cells of level $n$ over $S$. Clearly the cells are disjoint semi-algebraic sets. Therefore $\mathfrak{D}$ is indeed a cylindrical algebraic decomposition adapted to $\{V(I_k)\}$.
\end{proof}

It is time to introduce our cylindrical algebraic decomposition Algorithm based on the geometric tools we developed (see Algorithm \ref{algo-GeometricCAD}). The correctness of Algorithm \ref{algo-GeometricCAD} is obvious from Theorem \ref{thm-CAD-lifting}.

\begin{algorithm}[htbp]
	\KwIn{A list of varieties $\{V(I_k)\}_{k=1}^s$ in $\mathbb{A}_R^n$.}
	\KwOut{A cylindrical algebraic decomposition adapted to $\{V(I_1),\ldots,V(I_s)\}$.}
	\caption{GeometricCAD}\label{algo-GeometricCAD}
	
	\If(\tcp*[f]{The Projection Phase}){$n>1$}{
		$o_1,o_2:=\mathtt{CADProjection}(\{I_k\}_{k=1}^s)$\;
		$\mathfrak{D}:=\mathtt{GeometricCAD}(o_2)$ \tcp*{Compute a c.a.d.\ adapted to $o_2$}
	}
	\Else{		
		$\mathfrak{D}:=\{\bullet\}$ \tcp*{A c.a.d.\ of $R^0$}
	}
	$\mathtt{cells}:=\varnothing $ \tcp*{The cells of level $n$, initialized to be an empty list}
	\ForEach(\tcp*[f]{The Lifting Phase}){cell $S$ of level $n-1$ in $\mathfrak{D}$}{
		Choose a sample point $s\in S\subseteq R^{n-1}$\;
		$b:=\mathtt{False}$\;
		$\xi:=\varnothing$ \tcp*{The sections over $S$, initialized to be an empty list}
		\ForEach{ideal $I_k$}{
			$\ideal{f}$:=Specialize $I_k$ with $s\in S$ \tcp*{After specialization, the result is a principal ideal in $R[x_n]$}
			\If{$f\neq 0$}{
				$\{\rho_i\}:=$ The real roots of $f$ in ascending order\;
				\ForEach{$\rho_i$}{Append $\{(s,x)|s\in S,\allowbreak x \text{ is the i-th least real root of } I_k \text{ specialized at } s\}$ to $\xi$\;}
			}
			\Else(\tcp*[f]{$V(I_k)$ is vacuously geometrically delineable over $S$}){$b:=\mathtt{True}$\;}
		}
		Sort $\xi$ in ascending order, $\mathtt{cells}:=\mathtt{cells}\cup \xi$\;
		\If{$b$}{Append the bands over $S$ to $\mathtt{cells}$\;}		
	}
	Append $\mathtt{cells}$ to $\mathfrak{D}$,
	\Return{$\mathfrak{D}$}\;
	
	
\end{algorithm}

\section{Analysis}\label{sect-analysis}

\subsection{Connection to the Prior Work}

We compare our theory to some other existing methods. We hope the following discussion can be inspiring and illustrate the similarities and differences between our results and previous theories. 

\subsubsection{Collins's Original CAD}
It is very clear that our theory is deeply influenced by Collins's fundamental work \cite{collins1975quantifier}. Our algorithm follows the Projection-Lifting paradigm developed by Collins. Also, Algorithm \ref{algo-PROJ1} and \ref{algo-PROJ2} are analogues to Collins's Projection Operator \texttt{PROJ1} and \texttt{PROJ2}.

While Collins's goal is to compute an $\mathcal{F}$-sign-invariant cylindrical algebraic decomposition of the whole real affine space, our result is a c.a.d.\ adapted to a family of varieties. In order to compute an $\mathcal{F}$-sign-invariant\ c.a.d., Collins used his Projection Operator \texttt{PROJ} to reduce the task to computing a $\mathtt{PROJ}(\mathcal{F})$-sign-invariant c.a.d.\ in a lower-dimensional space. This idea is evolved in our geometric theory: a c.a.d.\ adapted to a family of varieties can be lifted from another lower-dimensional c.a.d.\ adapted to varieties. 

To exploit all equations, one has to study the geometry more carefully, so there is a significant difference in techniques involved. The mathematics behind Collins's work is the classical Elimination Theory based on resultants and subresultants, and our theory has more recent ingredients from modern algebraic geometry, which are translated into algorithms by Gr\"{o}bner basis via the standard algebra-geometry dictionary. These new ingredients allow us to manipulate varieties defined by ideals, rather than hypersurfaces defined by polynomials and reveal the geometric nature behind cylindrical algebraic decomposition. 

In fact, a sign-invariant c.a.d.\ can be thought as a c.a.d.\ adapted to the hypersurfaces (variety of a single polynomial). Therefore, Collins's Projection can be explained in our new theory as a special case that tackles hypersurfaces only. The Collins-Hong Projection Operator \texttt{PROJ} consists of three types of polynomials: leading coefficients, (sub-)discriminants and subresultants. By Corollary \ref{cor-collins-thm-1}, the leading coefficient of $f$ is used to construct the generic freeness of the map $\spec R[x_1,\ldots, x_{n+1}]/\ideal{f} \to \spec R[x_1,\ldots,x_n]$, the sub-discriminants of $f$ are used to count the geometric fiber of the same map. Similarly, the subresultants of $f$ and $g$ are used to construct the generic freeness of the map $\spec R[x_1,\ldots,x_{n+1}]/\ideal{f,g}\to \spec R[x_1,\ldots,x_n]$! So there are sufficiently many polynomials in the Collins-Hong Projection Operators to define the regions returned by our Algorithms \ref{algo-PROJ1} and \ref{algo-PROJ2}.

The advantage of the new idea is the ability to exploit all equations in the projection phase, avoiding unnecessarily refined decomposition. Actually, a $\mathtt{PROJ}(\mathcal{F})$-sign-invariant c.a.d.\ is not necessarily needed. Since we can conclude the geometric delineability over a region, it should be more economic if we work directly with regions instead of polynomials. The reason is simple: a $\mathtt{PROJ}(\mathcal{F})$-sign-invariant c.a.d.\ is adapted to all semi-algebraic sets that can be defined by $\mathtt{PROJ}(\mathcal{F})$, but a c.a.d.\ adapted to some semi-algebraic sets defined by $\mathtt{PROJ}(\mathcal{F})$ is not necessarily $\mathtt{PROJ}(\mathcal{F})$-sign-invariant. Therefore, a $\mathtt{PROJ}(\mathcal{F})$-sign-invariant c.a.d.\ might decompose the real space into more pieces than what is really needed. The following example we learned from Christopher Brown can be used to show the advantage of geometric theory.
\begin{example}\label{ex-5-var}
	Let $f=w^2 + ((x + y)^2 u - y z)\in R[x,y,z,u,w]$. We will build an $f$-sign-invariant cylindrical algebraic decomposition of $R^5$, using the classical method and our method.
	
	\begin{table}[h]
			\begin{center}
			\begin{tabularx}{0.87\textwidth}{c|c|c}
				\toprule
				& $\mathcal{F}$-sign-inv. c.a.d., where $\mathcal{F}=$ & c.a.d. adapted to \\
				\midrule
				Level 5 & $w^2+(x+y)^2 u - y z$ & $V(0), V(f)$  \\
				Level 4 & $(x+y)^2 u - y z$  & $V(0), V((x+y)^2 u - y z)$ \\
				Level 3 & $(x+y)^2,yz$ & $V(0), V(x+y), V(x+y,y z)$ \\
				Level 2 & $x+y,y$  & $V(0), V(x+y), V(x,y)$ \\
				Level 1 & $x$ & $V(0), V(x)$ \\
				\bottomrule
			\end{tabularx}
		\end{center}
		\caption{A comparison of computation of $f$-sign-invariant c.a.d.}\label{table-f-sign-inv.-cad}
	\end{table}
	
	Using the Collins-Hong Projection Operator \cite{collins1975quantifier,hong1990improvement}, the problem is reduced to computing lower-dimensional sign-invariant c.a.d.s of $R^4,\ldots,R^1$ and the projection polynomials can be found in the left column of Table \ref{table-f-sign-inv.-cad}. The final result consists of 297 cells in $R^5$.
	
	When applying McCallum-Brown Projection Operator \cite{mccallum1998improved,brown2001improved}, the same family of projection polynomials is generated. However, the nullification problem occurs, as the projection of $V((x+y)^2 u -y z)$ downto $\mathbb{A}_R^3$ is not quasi-finite (it has 1-dimensional fibers over $V(x+y,yz)$). In this case, the use of McCallum-Brown projection is not guaranteed to be correct by their theory (although it generates the same c.a.d.\ as Collins-Hong Projection does).
	
	Notice that a c.a.d.\ adapted to $\{V(f),V(0)\}$ is an $f$-sign-invariant c.a.d.\ of $R^5$, so it enables the use of our \texttt{GeometricCAD} algorithm. In our theory, this problem is reduced to compute a lower-dimensional c.a.d.\ adapted to another family of varieties recursively. The right column in Table \ref{table-f-sign-inv.-cad} shows the varieties needed to compute the c.a.d.. The final result consists of 75 cells in $R^5$. 
	
	Comparing 75 cells with 297 cells, we see a huge saving in cell numbers. Why does \texttt{GeometricCAD} generate a more compact c.a.d.? The reason lies in the geometric nature of cylindrical algebraic decomposition: sign invariance is more than necessary to conclude delineability (or geometric delineability). To see this, consider the c.a.d.\ constructed at level 3: a $\{(x+y)^2,y z\}$-sign-invariant c.a.d.\ $\mathfrak{D}_3$ versus a c.a.d.\ $\mathfrak{D}_3'$ adapted to $\{V(0),V(x+y),V(x+y,y z)\}$. Clearly $\mathfrak{D}_3$ is also adapted to these varieties. But $(x+y)^2 u - y z$ (or its locus) is obviously (geometrically) delineable over $V(0)\backslash V(x+y)$, $V(x+y)\backslash V(x+y, y z)$ or $V(x+y, y z)$ and the sign invariance of $\{(x+y)^2,yz\}$ is an overkill -- we do not really care what happens over $V(y z)$, as long as $x+y\ne 0$. So $\mathfrak{D}_3$ partitions $R^3$ into more cells than $\mathfrak{D}_3'$, resulting a larger decomposition in the end. As a consequence, the classical algorithm yields 3, 13, 39, 99 and 297 cells at each level, while our algorithm \texttt{GeometricCAD} generates 3, 9, 13, 25 and 75 cells at each level.
\end{example}

\subsubsection{CAD with Equational Constraints}

Usually, it is not necessary to decompose the whole affine real space but only a variety. For many problems from the real world, there are equational constraints that must be met.

Collins pointed out that the construction of cylindrical algebraic decomposition can be remarkably accelerated if taking equational constraints into consideration in \cite{collins1998quantifier}. Later McCallum made precise Collins's idea in \cite{mccallum1999projection} by reducing the projection polynomials in the first stage of projection phase. In \cite{mccallum2001propagation}, he further showed that repeated applications of ``semi-restricted projection'' in the projection phase is possible. The lower-dimensional equational constraints are obtained by a technique called ``EC Propagation''. To be more precise, Collins wrote in \cite{collins1998quantifier}:
\begin{quote}
	If $f_1$ and $f_2$ are both equational constraints polynomials then so is their resultant, $\res(f_1,f_2)$, since $$f_1=0\wedge f_2=0\implies\res(f_1,f_2)=0.$$
\end{quote}
And the propagation can be repeated many times. 

Our strategy to exploit equational constraints here is significantly different from Collins's idea. 
Since our theory is based on varieties, the ability to exploit equational constraints is intrinsic. If we need to decompose $V(I)$ into cells on which $g_1,\ldots,g_t$ are sign-invariant, then it suffices to compute a cylindrical algebraic decomposition adapted to 
$$\left\{V(I),V(I+\ideal{g_1}),\ldots,V(I+\ideal{g_t})\right\}.$$

We believe that our algorithm exploits equational constraints in a better way, as the following example demonstrates.

\begin{example}[{Example 7 in \cite{fukasaku2015real}}]\label{ex-FIS-7}
	The goal is to compute a cylindrical algebraic decomposition adapted to a one-dimensional variety (curve) $\mathcal{C}\subseteq \mathbb{A}_R^5=\spec R[x,y,z,u,w]$, and $\mathcal{C}$ is defined by the locus of polynomials 
	$$\begin{array}{rcl}
		f_1&=&w u-1,\\
		f_2&=&w^2-2 w y+u^2-2 u z-x, \\
		f_3&=&32 y^2+168 y z+40 y x-270 y+8 z^2+20 z x-390 z+4 x^2-105 x+450 \\
		f_4&=&240 y z-16 y x^2-532 y x-4480 y-800 z^3-1240 z^2 x-17720 z^2-408 z x^2\\
			& &-6214 z x+25240 z-40 x^3-550 x^2+5695 x+1050
		\text{ and }\\
		f_5&=&320 y z x+8320 y z+32 y x^2+264 y x-14840 y+240 z^2+16 z x^2-372 z x\\
		& &-23380 z-140 x^2 -2575 x+36750.
	\end{array}$$
	Our algorithm \texttt{GeometricCAD} finishes the computation in 10.67 seconds, generating 200 cells at level 5. 
	
	We tried to compute a c.a.d.\ adapted to $\mathcal{C}$ with the aid of QEPCAD \cite{brown2003qepcad}, which applies McCallum's equational constraint operator \cite{mccallum1999projection,mccallum2001propagation} after manually declaring $f_1,\ldots,f_5$ as equational constraints. However, QEPCAD quit after 14 minutes because of a software failure.
\end{example}

The use of resultants is completely abandoned in our geometric theory, and we believe this is beneficial because resultants do not always reflect the geometry faithfully. The phenomenon becomes apparent in EC Propagation as explained below.
\begin{enumerate}
	\item The locus of $\res_x(f,g)$ is not necessarily the closure of the projection of $V(f,g)$. In fact, we have 
	\begin{align*}\label{eqn-res-proj}V(\res_x(f,g))=V(\lc_x(f),\lc_x(g))\cup \pi(V(f,g))\\=V(\lc_x(f),\lc_x(g))\cup \overline{\pi(V(f,g))}\tag{$\star$},\end{align*}
	which follows from Proposition 6 and Corollary 7 in \cite[Chapter 3, Section 6]{cox2013ideals} quite easily. So there is a superfluous set $V(\lc_x(f),\lc_x(g))$, which does not necessarily lie in the silhouette of $V(f,g)$.
	\item Successively applying resultants may distort the geometry, making things even worse. Suppose $f,g,h$ are three polynomials in $x,y,z$. Set $r_1=\res_x(f,g)$, $r_2=\res_x(f,h)$ and $r=\res_y(r_1,r_2)$.
	Then $$\begin{array}{rcl}
		\pi_1(\pi_2(V(f,g,h)))&=&\pi_1(\pi_2(V(f,g)\cap V(f,h)))\\
		&\subseteq& \pi_1(\pi_2(V(f,g))\cap \pi_2(V(f,h)))\\
		&\subseteq& \pi_1(V(r_1)\cap V(r_2))\\
		&=& \pi_1(V(r_1,r_2))\\
		&\subseteq& V(r).
	\end{array}$$
	The first containment is by $\varphi(S\cap T)\subseteq \varphi(S)\cap \varphi(T)$ and the others follow from Equation (\ref{eqn-res-proj}). Each containment can be strict, which leads to an unnecessarily large estimate of the silhouette of $V(f,g,h)$ (usually reflected in superfluous factors of successive resultants).
\end{enumerate}

The failure of QEPCAD in Example \ref{ex-FIS-7} is possibly due to these reasons. The problem is overcome when using Gr\"{o}bner basis, because the vanishing set of the elimination ideal is always the closure of the projection.

%
%

\subsubsection{CAD by Regular Chain}
Apart from Collins's work, there is another approach to cylindrical algebraic decomposition, proposed by Changbo Chen, Marc Moreno Maza, Bican Xia and Lu Yang in \cite{chen2009computing} and \cite{chen2014incremental}, using Regular Chains. Our work was firstly motivated by \cite{chen2009computing} in the following way: since Regular Chains can be obtained from a lexicographic Gr\"{o}bner basis computation \cite{wang2020decomposition}, it is natural to ask if it is possible to build CAD using only Gr\"{o}bner Basis? 

Soon it turns out that \cite{chen2009computing} reflects some geometric nature of cylindrical algebraic decomposition and one of the highlights in \cite{chen2009computing} is that instead of directly building a cylindrical algebraic decomposition in $R^n$, they first build a cylindrical algebraic decomposition in $C^n$ then refine it to a cylindrical algebraic decomposition in $R^n$. From our perspective, a Regular Chain computation immediately provides generic freeness by inverting the leading coefficients and another discriminant computation counts the geometric fiber. 

But the dissimilarity is also clear. Thanks to the standard algebra-geometry dictionary, our theory can be directly translated into an algorithm, while the algorithm in \cite{chen2009computing} is sophisticated, using many high-level routines as black boxes due to the obscure relation between Regular Chains and varieties. Also, we can develop our theory much further without Cylindrical Projection, we do not use it until Theorem \ref{thm-freeness-implies-geometric-delineability}, which is about geometric delineability. It seems that recursively projecting the largest variable is an intrinsic nature of the theory of Regular Chains. This feature of Regular Chains makes it similar to the lexicographic Gr\"{o}bner Basis, which is usually hard to compute. As the readers can see, the monomial order does not play an important role in our work, and we can use cheaper monomial orders like graded reverse lexicographic order (grevlex) or block-grevlex (for elimination). Faug\`{e}re's \texttt{FGb} package supports these two monomial orders \cite{faugere2010fgb}, and is believed to be the state-of-the-art implementation of Gr\"{o}bner Basis.

\subsubsection{Preprocessing with Gr\"{o}bner Basis}
The idea of combining Gr\"{o}bner Basis and Cylindrical Algebraic Decomposition is actually not new. Buchberger and Hong has considered replacing equational constraints with their lexicographic Gr\"{o}bner Basis in \cite{buchberger1991speeding}, which usually speeds up the computation. It is also observed that, however sometimes computing a Gr\"{o}bner Basis would slow down the computation. Then, the experiments in \cite{buchberger1991speeding} was revisited by Wilson, Bradford and Davenport in \cite{wilson2012speeding}, with more advanced implementations of the algorithms. Huang \emph{et al.} applied machine learning to decide when to use a Gr\"{o}bner Basis in \cite{huang2016using}.

In short, the role that Gr\"{o}bner Basis plays in these works is a preprocessor, which is not the case in our theory. We use Gr\"{o}bner Basis as a fundamental algorithm to deal with ideals and varieties effectively.

\subsubsection{Comprehensive Gr\"{o}bner System} 
Suppose $I\unlhd k[\bs{y},\bs{x}]=k[y_1,\ldots,y_m,x_1,\ldots,x_n]$. A Comprehensive Gr\"{o}bner System for $I$ is a stratification of $\overline{k}^m$ into constructible sets $S_i$ such that each set is tagged with a finite subset $\mathcal{G}_i=\{f_{i,1},\ldots,f_{i,n_i}\}$ of $k[\bs{y},\bs{x}]$, satisfying for every $Y\in S_i$, the specialization of $\mathcal{G}_i$ at $Y$ is a Gr\"{o}bner Basis of the specialization of $I$ at $Y$. And a comprehensive Gr\"{o}bner Basis for $I$ is a finite subset of $k[\bs{y},\bs{x}]$ that is a Gr\"{o}bner Basis under all specializations. These concepts are proposed by Weispfenning in \cite{weispfenning1992comprehensive}. The theory of Comprehensive Gr\"{o}bner System/Basis is very rich and well suited to study parametric problems, readers may refer to \cite{lu2019survey} for a beautiful survey on this subject.

In \cite{weispfenning1998new}, Weispfenning used Comprehensive Gr\"{o}bner Bases and Hermite Quadratic Forms to develop a real quantifier elimination procedure. His approach roughly goes as:
\begin{enumerate}
	\item Convert the input to the prenex disjunctive normal form (DNF), which reduces the problem to producing an equivalent formula of $$(\exists x_1)\cdots(\exists x_k)(\bigwedge_i f_i(x_1,\ldots,x_n)=0\wedge\bigwedge_i h_i(x_1,\ldots,x_n)>0),$$
	which is the projection of a basic semi-algebraic set.
	\item Reduce to the generically finite case using Comprehensive Gr\"{o}bner Bases, then compute the Hermite Quadratic Forms associated with $1,h_i,h_i^2$.
	\item Then the number of real solutions can be read off from the sign variations of the characteristic polynomials of the quadratic forms (especially when the number is zero).
\end{enumerate}
This idea was further developed in \cite{fukasaku2015real}, where they used a more advanced Comprehensive Gr\"{o}bner System algorithm by Suzuki and Sato \cite{suzuki2006simple}.

We find that our Projection algorithms have a similar structure to the Comprehensive Gr\"{o}bner System construction in \cite{suzuki2006simple}: compute a reduced Gr\"{o}bner basis with respect to a block monomial order and pick out the leading coefficients. Moreover, we also use Hermite Quadratic Forms. These common features make our work comparable to \cite{weispfenning1998new,fukasaku2015real}. The major difference lies in the ways using Hermite Quadratic Forms. Weispfenning and his successors directly use the signature of the matrices to count the real roots, while we use the rank to count the complex roots (as Collins uses subresultants). We believe this can result in a significant difference in the number of basic semi-algebraic sets and the number of varieties, because the number of sign variance of the characteristic polynomial is far more difficult to describe than the rank (see \cite[Section 5.1]{fukasaku2015real}). Moreover, converting to disjunctive normal form will also produce many basic semi-algebraic sets, and many of them are actually empty. These problem are avoided in our theory because Cylindrical Algebraic Decomposition is able to postpone the construction of semi-algebraic sets in the Lifting Phase.

\subsection{Experimental Results}\label{subsect-experiments}

We report an implementation in this subsection and compare it to existing tools. The source codes and all examples in the experiment are available at 
\begin{center}
	\urlstyle{sf}
	\url{https://github.com/xiaxueqaq/GeometricCADv2}.
\end{center}
Our algorithm is implemented with \textsc{Mathematica} 13.3 \cite{Mathematica}. 
Our experiment runs on a Windows PC with Intel i7-11700 (2.50GHz) CPU and 32GB memory. 

We compare our implementation with \begin{enumerate}
	\item QEPCAD B 1.74 \cite{brown2003qepcad}, which applies McCallum-Brown Projection Operator \cite{mccallum1998improved,brown2001improved} by default. The command to launch QEPCAD is $$\texttt{qepcad +N500000000 +L2000000}.$$ Also, we declare all the Equational Constraints (whenever there are) in QEPCAD, so it can utilize the McCallum Equational Constraint Operator \cite{mccallum1999projection,mccallum2001propagation}.
	\item \texttt{CylindricalAlgebraicDecompose} \cite{chen2014cylindrical} in \textsc{Maple} 2021, which uses the algorithm in \cite{chen2014incremental} by Regular Chains. This tool is denoted by RC-CAD later.
\end{enumerate}

The testing examples are taken from \cite{bjorck1991faster}, \cite{england2020cylindrical}, \cite{fukasaku2015real} and \cite{hong1991comparison}, except 5-VAR, which is Example \ref{ex-5-var}. We summarize information about the benchmark in Table \ref{table-examples}, including the number of variables, polynomials, equations and the (Krull) dimension. The number of variables varies from 3 to 10 and the benchmark has both low-dimensional and high-dimensional varieties. Also, there are two families of scalable examples (Hong-$n$ and cyclic-$n$). Hence our benchmark is a comprehensive one.

The result is recorded in Table \ref{table-computation-time}. ``$>$1h'' stands for a killed computation that consumes more than 1 hour. ``Error'' means the computation halts due to a bug. Overall, it is evident that our algorithm solves the most problems. Also, our algorithm is of the best performance in 15 instances, which are 75\% of the benchmark.

In the first 6 examples, there is no equational constraints, and their goal is to build a sign-invariant c.a.d.\ to the input polynomials. The input of Hong-$n$ consists of two $n$-variate polynomials: $\sum_{i=1}^{n} x_i^2-1$ and $\prod_{i=1}^n x_i-1$. Although QEPCAD is the fastest tool on small instances like Hong-5 and Hong-6, it fails to finish the computation on larger problems. RC-CAD computes one more c.a.d.\ than QEPCAD, but it is slower than our algorithm. It crashes in Hong-8 for an unknown reason. Hence it is reasonable to say that our method outperforms the others in this category of data. While QEPCAD and GeometricCAD generate same numbers of cells on Hoon-$n$, it is interesting to see that RC-CAD uses fewer cells in c.a.d.\ construction.

Moving on to the other 14 examples, there are now equational constraints. The cyclic-$n$ family is a famous system in the symbolic computation community, which consists of cyclic sums of products of variables \cite{bjorck1991faster}. QEPCAD succeeds in 7 out of 14, but requires longer time in general. The failure of QEPCAD in FIS-7 and FIS-8 is due to an exhausted prime list in the \texttt{saclib} library. RC-CAD finishes two more instances and it is the most efficient tool on 3 examples, while our algorithm solves all the problems and uses the least time on 11 examples. Furthermore, neither QEPCAD nor RC-CAD is able to build a c.a.d.\ for cyclic-6, FIS-5, FIS-6 or FIS-7. As for the cell number, our algorithm always use less cells than QEPCAD, but RC-CAD's result is comparable to ours.


{\small
\begin{table}[htbp]
	\begin{center}
		\begin{tabular}{c|cccc}
			\toprule
			
			Examples & \#Variables & \#Polynomials & \#Equations & Dimension \\
			
			\midrule
			
			5-VAR & 5 & 1 & 0 & 5 \\			
			
			Hong-4 & 4 & 2 & 0 & 4 \\
			
			Hong-5 & 5 & 2 & 0 & 5 \\
			
			Hong-6 & 6 & 2 & 0 & 6 \\
			
			Hong-7 & 7 & 2 & 0 & 7 \\
			
			Hong-8 & 8 & 2 & 0 & 8 \\
			
%
%
%

			cyclic-4 & 4 & 4 & 4 & 0 \\
			
			cyclic-5 & 5 & 5 & 5 & 0 \\
			
			cyclic-6 & 6 & 6 & 6 & 0 \\			
			
			EBD-2 & 3 & 3 & 2 & 1 \\
			
			EBD-5 & 5 & 6 & 4 & 1 \\
			
			EBD-8 & 4 & 4 & 3 & 1 \\
			
			FIS-1 & 5 & 2 & 1 & 4 \\
			
			FIS-2 & 4 & 3 & 3 & 1 \\
			
			FIS-3 & 7 & 4 & 4 & 3 \\
			
			FIS-4 & 3 & 3 & 2 & 1 \\
			
			FIS-5 & 5 & 5 & 5 & 0 \\
			
			FIS-6 & 5 & 3 & 3 & 2 \\
			
			FIS-7 & 5 & 6 & 5 & 1 \\
			
			FIS-8 & 10 & 15 & 12 & 0 \\
			
			\bottomrule
	\end{tabular}			\end{center}
	\caption{Information about the Benchmark}\label{table-examples}
\end{table}

\begin{table}[htbp]
	\begin{center}
		\begin{tabular}{c|cc|cc|cc}
			\toprule
			Tools& \multicolumn{2}{c|}{QEPCAD B} & \multicolumn{2}{c|}{RC-CAD}&  \multicolumn{2}{c}{GeometricCAD}\\
			& Time (s) & \#Cells & Time (s) & \#Cells  & Time (s) & \#Cells \\
			\midrule	
			
			5-VAR & 1.50 & 297 & 0.09 & 75 & \textbf{0.08} & 75 \\					
			
			Hong-4 & 1.35 & 575 &  1.19 & 491  & \textbf{0.56} & 575\\
			
			Hong-5 & \textbf{1.35} & 2241 &  5.16  & 1731  & 2.53 & 2241\\
			
			Hong-6 & \textbf{1.70} & 7967 & 20.00  & 5631   & 9.77 & 7967\\
			
			Hong-7 & $>$1h & N/A &  133.86 & 19743   & \textbf{49.52} & 29841\\
			
			Hong-8 & $>$1h & N/A   & Error  & N/A  & \textbf{270.33} & 103343 \\
			
%
%
%

			cyclic-4 & 1.56 & 129 & 0.06 & 4 & \textbf{0.03} & 4 \\
			
			cyclic-5 & $>$1h & N/A & 4.58 & 10 & \textbf{0.58} & 10 \\
			
			cyclic-6 & $>$1h & N/A & $>$1h & N/A & \textbf{5.48} & 24 \\
			
			EBD-2 & 1.17 & 41 &  \textbf{0.02}  & 4  & \textbf{0.02} & 5\\
			
			EBD-5 & 1.33 & 107 &  0.27 & 13   & \textbf{0.05} & 10 \\
			
			EBD-8 & 107.89 & 257 & $>$1h  & N/A & \textbf{0.23} & 4 \\										
			
			FIS-1 & \textbf{3.68} & 114541 &  40.25 & 5333 & 16.11 & 4063\\
			
			FIS-2 & 1101.60 & 5467 &  \textbf{1.59} & 79 & 2.66 & 299 \\
			
			FIS-3 & $>$1h & N/A &  \textbf{644.91} & 16928 & 1654.84 & 66704\\
			
			FIS-4 & 4.76 & 227 &  1.83 & 24 & \textbf{0.19} & 16\\
			
			FIS-5 & $>$1h & N/A &  $>$1h & N/A & \textbf{0.20} & 3 \\
			
			FIS-6 & $>$1h & N/A &  $>$1h & N/A & \textbf{2834.00}  & 32875 \\
			
			FIS-7 & Error & N/A &  $>$1h & N/A &  \textbf{10.67} & 200 \\
			
			FIS-8 & Error & N/A & 122.58 & 3 & \textbf{12.75} & 3\\
			\bottomrule
			
		\end{tabular}			
	\end{center}
	\caption{Computation Time and Number of Cells}\label{table-computation-time}
\end{table}}

We have not included any complexity analysis because,
\begin{enumerate}
	\item The decision problem over the real closed field is in \EXPSPACE~and it is also conjectured that this is \EXPSPACE-complete \cite{ben1984complexity}. Moreover, it is shown in \cite{davenport1988real,weispfenning1988complexity,brown2007complexity} that a cylindrical algebraic decomposition needs doubly-exponential number of cells in the worst case. So our algorithm probably does not improve the worst-case doubly-exponential bound of a cylindrical algebraic decomposition construction.
	\item Our algorithm heavily relies on Gr\"{o}bner Basis. But the complexity of Buchberger's Algorithm is still not well understood in the geometric context. While the ideal membership problem is shown to be \EXPSPACE-complete \cite{MAYR1982complexity,Mayr1989membership}, Buchberger's algorithm and its variants are often able to tackle problems of large size. In fact, there are many works to study why Gr\"{o}bner basis computation is feasible in the cases of real interest in algebraic geometry, and sharper upper bounds are obtained by analyzing zero-dimensional systems \cite{Lakshman1991A, Lakshman1991on}, low-dimensional systems \cite{mayr2010degree,mayr2013dimension}, regular sequences (Faug\`{e}re's F5) \cite{faugere2002new, BARDET2015on}.
	
\end{enumerate}

Still, the experimental data shows that our method goes beyond the reach of classical algorithms.

\subsection{Future Work}

There are many possible directions for future researches. We list some of them below.

\subsubsection{Adapting Improved Projection Operators to the Geometric Context}

There are many improved projection operators besides Hong's \cite{hong1990improvement}. McCallum proposed his simplified projection operator \cite{mccallum1988improved,mccallum1998improved} and Brown further improved McCallum's work \cite{brown2001improved}, although their lifting could fail in rare cases (the nullification problem). Lazard came up with another kind of projection operator and he used a slightly different lifting scheme \cite{lazard1994improved}. But he did not provide a proof for an essential claim in his paper. Lazard's CAD construction was only validated recently in \cite{mccallum2019validity}.

Following this paper, it is natural to ask, what is the geometric theory behind McCallum-Brown Projection and Lazard Projection? How can they reduce the Projection Operator size? Can we adapt their work to the Geometric CAD? The answers to these problems will surely help researchers gain a deeper understanding of the geometry of CAD and design more efficient algorithms.

%
%

\subsubsection{Full Quantifier Elimination}

A Quantifier Elimination (QE) algorithm is an algorithm that takes a quantified sentence and returns a quantifier-free sentence that is equivalent to the input. Collins's method is the first practical Quantifier Elimination algorithm for $\mathbb{R}$. It performs an additional step after the construction of CAD to find an equivalent quantifier-free formula.

Strictly speaking, our algorithm here does not directly eliminate quantifiers, because our cells are described by semi-algebraic continuous functions, which are generally not polynomials. However, it does have the ability of describing the regions in an extended language (so called ``Extended Tarski Formula''). In other words, if we allow ``root functions'' in the language, then our algorithm can give a quantifier elimination procedure.

Collins's original CAD can be refined to give a polynomial-only description of the cells, using more polynomials ($\mathtt{APROJ}$, the augmented projection in \cite{collins1975quantifier}) and Thom's Encoding. A possible line of research is to develop a Geometric CAD theory of quantifier elimination.

\section*{Acknowledgment}
 The author would like to thank Changbo Chen, James Davenport, Tereso del Río Almajano, Matthew England, Jingjun Han, Scott McCallum, Mohab Safey El Din and Longke Tang for inspiring discussions. Also, the author would like to thank Manuel Kauers for generously allowing us to modify and redistribute his package in an earlier implementation. The author would like to thank Christopher Brown for instruction on using QEPCAD. The author would like to thank Hoon Hong for spending his spare time helping with revision. The author would like to thank the referees for their valuable comments and suggestions. Finally, the author would like to express his heartfelt appreciation to his advisor, Bican Xia, for continuous guidance and support. Without him, this paper would not exist.

\bibliographystyle{alpha}
\bibliography{bibs}

\appendix

\end{document}